\documentclass[english, a4paper, 10pt]{article}

\usepackage[english]{babel}
\usepackage[T1]{fontenc}
\usepackage{authblk}

\usepackage[a4paper,top=3cm,bottom=3cm,left=3cm,right=3cm,marginparwidth=1.75cm]{geometry}

\usepackage{empheq}
\usepackage{graphicx}
\usepackage[colorinlistoftodos]{todonotes}
\usepackage{amsmath}
\usepackage{mathtools}
\usepackage{upgreek}
\usepackage{amssymb}
\usepackage{amsfonts}
\usepackage{amsthm}
\usepackage{enumitem}
\usepackage{mathrsfs}
\usepackage{fontenc}
\usepackage{verbatim}
\usepackage{framed}
\usepackage{algorithm}
\usepackage{bbm}
\usepackage{algpseudocode}
\usepackage{appendix}
\usepackage{color}
\usepackage{multirow}
\usepackage[colorlinks=true, allcolors=blue]{hyperref}
\usepackage{caption}
\usepackage{mfirstuc}
\usepackage[figuresright]{rotating}

\numberwithin{equation}{section}
\theoremstyle{plain}
\newtheorem{thm}{Theorem}

\newtheorem{lem}{Lemma}

\newtheorem{prom}{Problem}

\newtheorem{assu}{Assumption}[section]

\newtheorem*{rem}{Remark}

\theoremstyle{definition}

\theoremstyle{remark}

\title{A deep learning method for solving stochastic optimal control problems driven by fully-coupled FBSDEs}

\author{
    Shaolin Ji$^1$, Shige Peng$^2$, Ying Peng$^1$, and Xichuan Zhang$^3$
    \bigskip
    \\
    \small{$^1$Zhongtai Securities Institute for Financial Studies, Shandong University, 250100, China}
    \\
    E-mail: \{jsl, pengy\}@sdu.edu.cn
    \smallskip
    \\
    \small{$^2$School of Mathematics, Shandong University, 250100, China}
    \\
    E-mail: peng@sdu.edu.cn
    \\
    \small{$^3$Intelligent Game and Decision Lab (IGDL), 100091, China}
    \\
    E-mail: zhxc@alu.hit.edu.cn
}

\begin{document}

\maketitle

\begin{abstract}

In this paper,
we mainly focus on the numerical solution of high-dimensional stochastic optimal control problem driven by fully-coupled forward-backward stochastic differential equations (FBSDEs in short) through deep learning. We first transform the problem into a stochastic Stackelberg differential game problem (leader-follower problem), then a bi-level optimization method is developed where the leader's cost functional and the follower's cost functional are optimized alternatively via deep neural networks. As for the numerical results, we compute two examples of the investment-consumption problem solved through stochastic recursive utility models,
and the results of both examples demonstrate the effectiveness of our proposed algorithm.

\end{abstract}


\section{Introduction}\label{sec-intro}

Bismut \cite{bismut1973conjugate} first introduced linear backward stochastic differential equations (BSDEs in short) as the adjoint equation of the classical stochastic optimal control problem.
In 1990, Pardoux and Peng firstly proved the existence and uniqueness of nonlinear BSDEs with Lipschitz condition \cite{pardoux1990adapted}.
Since then,
the theory of BSDEs has been studied by many researchers and applied in a wide range of areas,
such as in stochastic optimal control and mathematical finance~\cite{peng1997financial,peng1990stochasticmax}.
When a BSDE is coupled with a (forward) stochastic differential equation (SDE in short),
the system is usually called a forward-backward stochastic differential equation (FBSDE in short).
We can refer to the literatures in ~\cite{wu1999Fully,antonelli1993backward,ma1996hedging,hu1995solution,pardoux1999forward,ma1994solving} which studied the existence, uniqueness and the applications of coupled or fully-coupled FBSDEs.

In 1993, Peng \cite{peng1993backward} first established a local stochastic maximum
principle (SMP) for stochastic optimal control problems driven by FBSDEs. Then, the
local SMP for other related problems were studied by
Dokuchaev and Zhou \cite{dokuchaev1999JMAA}, Ji and Zhou \cite{ji2006cis}, Yong and Zhou \cite{Yong_stochastic_control} (
see also the references therein). When the control domain is non-convex, Hu  \cite{hu2017puqr}, Hu et al. \cite{hu2018global} built the global stochastic maximum
principle for stochastic optimal control problems driven by BSDEs and FBSDEs respectively.
On the other hand, the dynamic programming principle (DPP) and related
Hamilton-Jacobi-Bellman (HJB) equations have been intensively studied by Li and Wei \cite{li2014siam-jco}, Hu et al.  \cite{hu2019siamCO} for
this kind of stochastic optimal control problems. Furthermore, Hu et al. \cite{hu2020esaim-cocv} revealed the relationship between the SMP and the DPP for a
stochastic optimal control problem where the system is governed by a
 fully coupled FBSDE.
However, few literature has studied the numerical method for solving the stochastic optimal control problems driven by FBSDEs.


In this paper,
we aim to solve the high-dimensional stochastic optimal control problem where the state equation is governed by
\begin{equation}\label{FBControlSys-1}
    \left\{
    \begin{array}{l}
        d x(t) = b(t,x(t),y(t),z(t),u(t)) d t + \sigma(t,x(t),y(t),z(t),u(t)) d B(t), \vspace{1ex} \\
        -d y(t) = l(t,x(t),y(t),z(t),u(t)) d t - z(t) d B(t), \vspace{1ex} \\
        x(0) = a, \qquad y(T) = g(x(T)).
    \end{array}
    \right.
\end{equation}
and the cost functional is defined by
\begin{equation}\label{CostFB-1}
    J(u(\cdot)) = E \left[\int_0^T f(t,x(t),y(t),z(t),u(t)) d t + h(x(T)) + \gamma(y(0))\right] .
\end{equation}

Recall that for the high-dimensional classical stochastic optimal control problem where the state equation is a SDE, Han and E \cite{han2016deep} solved it through deep learning method. In more details, they
developed a feed-forward neural network to approximate the control $u(\cdot)$ and regard the cost functional $J(u(\cdot))$ as the optimization objective (see also \cite{carmona2021convergence2}). Comparing with the classical stochastic optimal control problem, the difficulty of our problem lies in that we need to deal with the high-dimensional FBSDE \eqref{FBControlSys-1} firstly. And it is well-known that to obtain the numerical solution of \eqref{FBControlSys-1} itself is a difficult problem.

The traditional methods for solving the FBSDEs include the partial differential equation (PDE in short) methods and the probabilistic methods, such as \cite{numericalPDE2012review,regression_for_BSDEs,FBSDE_Ma2008,time_discretize_FBSDE_Zhang,Zhao2016Multistep,Milstein2004Numerical,Yu2016Efficient,huijskens2016efficient}.
However, most of these methods can not deal with high-dimensional FBSDEs.
In this paper, we reformulate
the fully-coupled FBSDE \eqref{FBControlSys-1} as a stochastic optimal control problem with the cost functional given by the error between the terminal condition and the solution of  the FBSDE,
$$J(y_0,z(\cdot)) = E \left[|y(T)-g(x(T))|^2\right] ,$$
and build two feed-forward neural networks to approximate the control variables $(y_0,z(\cdot))$.
This idea for solving \eqref{FBControlSys-1} is currently a common method in solving high-dimensional BSDEs and FBSDEs (see e.g.~\cite{WeinanDLforHigh_dim,WeinanDLforBSDE,deeplearning_FBSDE,Peng_FBSDE_numerical,hure_deep_2020,pham_neural_2021,beck_machine_2019}).

Based on the above analysis, we put forward the following ideas to solve our stochastic optimal control problem \eqref{FBControlSys-1}-\eqref{CostFB-1} by a novel deep learning method. Firstly, we transform \eqref{FBControlSys-1}-\eqref{CostFB-1} into a stochastic Stackelberg differential game problem (leader-follower problem).
The goal of the follower is to find a pair of optimal control $(y_0,z(\cdot))$ which minimizes the functional $J(y_0,z(\cdot))$ under a given control $u(\cdot)$,
while the goal of the leader is to find an optimal control $u(\cdot)$ which minimizes the cost functional $J(u(\cdot))$.
Secondly, a bi-level optimization method is introduced to solve the high-dimensional stochastic Stackelberg differential game problem.

Specifically, we construct a feed-forward neural network to approximate the control $u(\cdot)$
and regard the cost functional $J(u(\cdot))$ as the optimization objective for the leader.
We also build two feed-forward neural networks to approximate the control $(y_0,z(\cdot))$ and regard $J(y_0,z(\cdot))$ as the optimization objective for the follower.
As for the parameter updates of the neural networks, we update the parameters approximating $(y_0,z(\cdot))$ through $J(y_0,z(\cdot))$ and the parameters approximating $u(\cdot)$ through $J(u(\cdot))$ alternatively. In more details, in each training step, we first fix an approximated $u(\cdot)$ to perform $\kappa$ $(\kappa\geq 1,\kappa\in\mathbb{N}^+)$ times updates for the parameters approximating $(y_0,z(\cdot))$, then based on the obtained $(y_0,z(\cdot))$, we perform  one time update for the parameters approximating $u(\cdot)$. The training will be repeated until a convergence result is obtained.
Note that under a given $u(\cdot)$, the follower should choose the optimal $(y_0,z(\cdot))$ such that $J(y_0,z(\cdot))=0$.
Here in order to improve the computation efficiency, we use a relaxing method in the parameter update process.
That is, instead of forcing $J(y_0,z(\cdot))=0$, we control the number of updates by choosing an enough large coefficient $\kappa$ so that $J(y_0,z(\cdot))$ can obtain a sufficiently small value which is close to $0$, but not strictly equal to $0$.

To show the feasibility of our proposed method,
we compute the investment-consumption problem for stochastic recursive utilities in a financial market. The concept of stochastic recursive utility was first introduced by Duffie and Epstein \cite{duffie1992stochastic}.
In fact, the stochastic recursive utility is associated with the solution of a
particular backward stochastic differential equation (BSDE).
From the BSDE point of view, El Karoui et al.
\cite{peng1997financial,el2001dynamic} considered a more general class of recursive utilities defined as
the solutions of BSDEs.
In this paper, we compute two examples. In the first example, the stochastic recursive utility is described by a linear BSDE. For this case, the investment-consumption problem can be transformed to a classical stochastic optimal control problem.
In this way, we can compare our computation results with those obtained by the numerical methods (such as the method in \cite{han2016deep}) for solving classical stochastic optimal control problems.
The numerical results show that the value of the recursive utility functions obtained by the two methods are very close.
In the second example,
we consider a more general recursive utility problem whose generator contains the $z$ term (refer to \cite{peng1997financial,el2001dynamic,chen2002ambiguity}).
Chen and Epstein \cite{chen2002ambiguity} studied the recursive utility problem with $z$ term,
in which $z$ is regarded as the volatility of utility.
The computation results show that our proposed method is also effective for this general case.





This paper is organized as follows.
In \autoref{sec-form}, we describe the stochastic optimal control problem driven by fully-coupled FBSDEs and reformulate it as a stochastic Stackelberg differential game problem. In \autoref{sec-schemes}, we present our proposed bi-level optimization method for solving the stochastic Stackelberg differential game problem via deep learning. As the numerical results, we compute the investment-consumption problem for stochastic recursive utilities in \autoref{sec-numerical-results}.
In \autoref{sec:conclusion}, we make a brief conclusion.

\section{Statement of the problem}\label{sec-form}

\subsection{The stochastic control problem driven by FBSDE}\label{ssec-Prob}

On a given complete probability space $(\Omega, \mathcal{F}, \mathbb{P})$,
let $(B(t), t\geq 0)$ be a standard $d$-dimensional Brownian motion defined on a finite interval $[0,T]$ for some given constant $T>0$.
We denote
\begin{equation*}
     \mathcal{F}_t = \sigma\{B(s), 0\leq s\leq t\}
 \end{equation*}
the natural filtration of $\{B(t)\}_{0\leq t\leq T}$ and $\mathbb{F} := \{\mathcal{F}_t\}_{0\leq t\leq T}$.
$\mathcal{F}_0$ contains all the $\mathbb{P}$-null set in $\mathcal{F}$ and $\{\mathcal{F}_t\}_{t\ge 0}$ is right continuous.
$M^2(0,T;\mathbb{R}^n)$ is denoted as the space of all mean square-integrable $\mathcal{F}_t$-adapted and $\mathbb{R}^n$-valued processes.
It is a Hilbert space with the norm
\begin{equation*}
    \|v(\cdot)\|_{M^2(0,T;\mathbb{R}^n)} = \Big( \mathbb{E} \big[ \int_0^T |v(t)|^2 d t \big] \Big)^{1/2}.
\end{equation*}
Set
\begin{equation*}
L^2(\Omega, \mathcal{F}_t, \mathbb{P}) = \left\{ \xi|\xi \in \mathbb{R}^n \mbox{ is } \mathcal{F}_t\mbox{-measurable and } \mathbb{E}\left[|\xi|^2\right] <\infty\right\},
\end{equation*}
where $|x|$ is denoted as the Euclidean norm of $x\in\mathbb{R}^n$ and $\left\langle x,y \right\rangle$ is denoted as the Euclidean inner product of $x,y\in\mathbb{R}^n$.



In this paper,
we consider the following controlled fully-coupled FBSDE,
\begin{equation}\label{FBControlSys}
    \left\{
    \begin{array}{l}
        d x(t) = b(t,x(t),y(t),z(t),u(t)) d t + \sigma(t,x(t),y(t),z(t),u(t)) d B(t), \vspace{1ex} \\
        -d y(t) = l(t,x(t),y(t),z(t),u(t)) d t - z(t) d B(t), \vspace{1ex} \\
        x(0) = a, \qquad y(T) = g(x(T)),
    \end{array}
    \right.
\end{equation}
where
\begin{align*}
    b: & \quad [0,T] \times \mathbb{R}^{n} \times \mathbb{R}^{m} \times \mathbb{R}^{m\times d} \times U \rightarrow \mathbb{R}^{n}, \\
    \sigma: & \quad [0,T] \times \mathbb{R}^{n} \times \mathbb{R}^{m} \times \mathbb{R}^{m\times d} \times U \rightarrow \mathbb{R}^{n\times d}, \\
    l: & \quad [0,T] \times \mathbb{R}^{n} \times \mathbb{R}^{m} \times \mathbb{R}^{m\times d} \times U \rightarrow \mathbb{R}^{m}, \\
    g: & \quad \mathbb{R}^n \rightarrow \mathbb{R}^m,
\end{align*}
are given $C^1$ functions with respect to $(x,y,z,u)$. The process $u(\cdot)$ taking value in a given nonempty convex set $U \subseteq \mathbb{R}^{k} $ in the system \eqref{FBControlSys} is called an admissible control.
For a given admissible control $u(\cdot)$,
the corresponding solution of the system \eqref{FBControlSys} is denoted by
\[
  (x(t),y(t),z(t)) = (x^{u}(t),y^{u}(t),z^{u}(t)), \quad \mbox{ for all} \quad t\in[0,T].
\]
The set of all admissible control $u(\cdot)$ is denoted as $\mathcal{U}_{ad}[0,T]$,
and the corresponding 4-tuple $(x(t),y(t),z(t),u(t))$ is called an admissible 4-tuple.

For a given control $u(\cdot)$,
$b,\sigma,l,g$ are functions of $(t,x,y,z)$.
For these functions, we need to introduce some assumptions.

Given an $m\times n$ full-rank matrix $G$ and set
\begin{equation*}
    w =
    \begin{pmatrix}
        x\\
        y\\
        z
    \end{pmatrix} \in \mathbb{R}^n \times \mathbb{R}^m \times \mathbb{R}^{m\times d}, \qquad A(t,w) =
    \begin{pmatrix}
        -G^Tl \\
        Gb\\
        G\sigma
    \end{pmatrix}(t,w),
\end{equation*}
where $G\sigma = (G\sigma_1,G\sigma_2,\cdots,G\sigma_d)$, and
\begin{equation*}
    \left\langle w^1, w^2 \right\rangle = \left\langle x^1, x^2 \right\rangle + \left\langle y^1, y^2 \right\rangle + \left\langle z^1, z^2 \right\rangle.
\end{equation*}

\begin{assu}\label{assu-1}
  \begin{itemize}
    \item[(i)] $A(t,w)$ is uniformly Lipschitz with respect to $w$;
    \item[(ii)] for each $w \in \mathbb{R}^n \times \mathbb{R}^m \times \mathbb{R}^{m\times d}$, $A(\cdot,w) \in M^2(0,T;\mathbb{R}^n)$;
    \item[(iii)] $g(x)$ is uniformly Lipschitz with respect to $x \in \mathbb{R}^n$.
  \end{itemize}
\end{assu}
The following monotonic conditions firstly introduced in \cite{wu1999Fully} are also needed.
\begin{assu}\label{assu-2}
  \begin{equation}
    \begin{aligned}
      \left\langle A(t,w)- A(t,w'), w-w' \right\rangle & \leq -\beta_1 |G\hat{x}|^2 - \beta_2(|G^T\hat{y}|^2+|G^T\hat{z}|^2), \\
      \left\langle g(x) - g(x'), G(x-x') \right\rangle & \geq \mu_1 |G\hat{x}|^2, \\
      \forall w = (x,y,z), \ w'=(x',y',z'), & \ \hat{x} = x-x', \ \hat{y} = y-y', \ \hat{z} = z-z',
    \end{aligned}
  \end{equation}
  or
  \begin{equation}
    \begin{aligned}
      \left\langle A(t,w)- A(t,w'), w-w' \right\rangle & \geq \beta_1 |G\hat{x}|^2 + \beta_2(|G^T\hat{y}|^2+|G^T\hat{z}|^2), \\
      \left\langle g(x) - g(x'), G(x-x') \right\rangle & \leq -\mu_1 |G\hat{x}|^2, \\
      \forall w = (x,y,z), \ w'=(x',y',z'), & \ \hat{x} = x-x', \ \hat{y} = y-y', \ \hat{z} = z-z',
    \end{aligned}
  \end{equation}
  where $\beta_1,\beta_2$ and $\mu_1$ are non-negative constants with $\beta_1 + \beta_2 >0$, $\beta_2+\mu_1>0$.
  Moreover we have $\beta_1>0,\mu_1>0$ (resp. $ \beta_2>0 $) when $m>n$ (resp. $m<n$).
\end{assu}

We also introduce the following Lemma \ref{lem-1} without proof,
and the proof of this lemma can be found in \cite{wu1999Fully}.

\begin{lem}\label{lem-1}
  For any given admissible control $u(\cdot)$,
  let \autoref{assu-1} and \autoref{assu-2} hold.
  Then the FBSDE \eqref{FBControlSys} has a unique adapted solution $(x(\cdot),y(\cdot),z(\cdot))$.
\end{lem}

In the following, we define the cost functional of the stochastic optimal control problem by
\begin{equation}\label{CostFB}
    J(u(\cdot)) := E \left[\int_0^T f(t,x(t),y(t),z(t),u(t)) d t + h(x(T)) + \gamma(y(0))\right] ,
\end{equation}
where
\begin{align*}
    f: & \quad [0,T] \times \mathbb{R}^{n} \times \mathbb{R}^{m} \times \mathbb{R}^{m\times d} \times U \rightarrow \mathbb{R}, \\
    h: & \quad \mathbb{R}^{n} \rightarrow \mathbb{R}, \\
    \gamma: & \quad \mathbb{R}^{m} \rightarrow \mathbb{R},
\end{align*}
are given $C^1$ functions with respect to $(x,y,z,u)$.
The 4-tuple process $(x(t),y(t),z(t),u(t))$ satisfies equation \eqref{FBControlSys}.
Then the stochastic control problem driven by FBSDE can be described as following.

\begin{prom}\label{prom-1}
  The optimal control problem is to find an admissible control $u^*(\cdot)$ over $\mathcal{U}_{ad}[0,T]$ such that
  \begin{equation}\label{promFBSC}
      J(u^*(\cdot)) = \inf_{u(\cdot)\in \mathcal{U}_{ad}[0,T]} J(u(\cdot)).
  \end{equation}
\end{prom}

If there exists a control $u^*(\cdot)$ which minimizes $J(u(\cdot))$ over $\mathcal{U}_{ad}[0,T]$,
 then we call it an optimal control. \eqref{FBControlSys} is called the optimal state equation and its solution $(x^*(\cdot), y^*(\cdot), z^*(\cdot))$ is called an optimal trajectory. As the existence of stochastic control problems driven by FBSDEs is a difficult problem (\cite{buckdahn2010existence},\cite{bahlali2011existence},\cite{bahlali2018existence}), and the relaxing method should be used even for some simple stochastic control problems such as the Linear Quadric problem. Therefore, in order to focus on the numerical method for solving \autoref{prom-1}, here we assume that the solution of \autoref{prom-1} exists.

\begin{rem}
    In Problem \ref{prom-1} and Lemma \ref{lem-1},
    we study the stochastic control problem driven by fully-coupled FBSDE \eqref{FBControlSys},
    which needs some strong conditions (such as the monotonic conditions \autoref{assu-1} and \autoref{assu-2}).
    If the forward SDE in \eqref{FBControlSys} does not contain $y$ and $z$ terms,
    such strong conditions can be relaxed \cite{peng1993backward}.
\end{rem}

\subsection{The Equivalent problem}\label{ssec-Equal}

As discussed in \cite{KohlmannZhou2000,LimZhou2001},
$z$ can be regarded as a control variable and $y(T)=g(x(T))$ in \eqref{FBControlSys} can be seen as a terminal state constraint.
Then \autoref{prom-1} can be reformulated as the following stochastic control problem to minimize
\begin{equation}\label{Wu-Yong}
    J(u(\cdot),y_0,z(\cdot)) = E \left[\int_0^T f(t,x(t),y(t),z(t),u(t)) d t + h(x(T)) + \gamma(y(0))\right] ,
\end{equation}
over $u(\cdot)\in\mathcal{U}_{ad}[0,T],y_0\in\mathbb{R}^m,z(\cdot)\in M^2(0,T;\mathbb{R}^{m\times d})$,
the state equation is given as
\begin{equation}\label{forFBSCsys}
    \left\{
    \begin{array}{l}
      dx(t) = b(t,x(t),y(t),z(t),u(t)) d t + \sigma(t,x(t),y(t),z(t),u(t)) d B(t), \vspace{1ex} \\
      -dy(t)= l(t,x(t),y(t),z(t),u(t)) d t - z(t) d B(t), \vspace{1ex} \\
      x(0) = a, \qquad y(0) = y_0,
    \end{array}
    \right.
\end{equation}
with the terminal state constraint
\begin{equation}\label{state-cons}
    E\left[|y(T)-g(x(T))|^2\right]=0.
\end{equation}

In this paper, our idea is to relax the terminal state constraint as an optimization objective which must be met firstly. Following this idea, the above stochastic control problem \eqref{Wu-Yong} becomes a stochastic Stackelberg differential game problem (leader-follower problem). In more details, we
assume that there is one follower and one leader, and their cost functionals are given as follows
\begin{equation}\label{CostFB-game}
    \begin{aligned}
        J_1(u(\cdot),y_0,z(\cdot)) &:= E \left[|y(T)-g(x(T))|^2\right],\\
        J_2(u(\cdot),y_0,z(\cdot)) &:= E \left[\int_0^T f(t,x(t),y(t),z(t),u(t)) d t + h(x(T)) + \gamma(y(0))\right] .
    \end{aligned}
\end{equation}
For any choice $u(\cdot)$ of the leader and a fixed initial state $a$, the goal of the follower is to minimize the functional $J_1(u(\cdot),y_0,z(\cdot))$ over $(y_0,z(\cdot)) \in \mathcal{A}$ where $\mathcal{A}$ is the set of all admissible $(y_0,z(\cdot))$. Denote the optimal $(y_0,z(\cdot))$ of the follower as $(y_0^u,z^u(\cdot))$ which clearly depends on $u(\cdot)$ of the leader. Then, the goal of the leader is to minimize the cost functional $J_2(u(\cdot),y_0^u,z^u(\cdot))$ over $u(\cdot)\in \mathcal{U}_{ad}[0,T]$. More details for describing the Stackelberg differential game problem can be found in \cite{yong2002leader}, with the state equation being a linear Ito-type stochastic differential equation and the cost functionals being quadratic.

In a more rigorous way, for $b,\sigma,l,g$ in \eqref{forFBSCsys} and $f,h,\gamma$ in \eqref{CostFB-game},
we introduce the following additional assumptions.

\begin{assu}\label{assu-3}
  \begin{itemize}
    \item[(i)] $b,\sigma,l,g,f,h$ and $\gamma$ are continuously differential;
    \item[(ii)] the derivatives of $b, \sigma, l, g$ are bounded;
    \item[(iii)] the derivatives of $f$ are bounded by $C(1+|x|+|y|+|z|+|u|)$;
    \item[(iv)] the derivatives of $h$ and $\gamma$ are bounded by $C(1+|x|)$ and $C(1+|y|)$, respectively.
  \end{itemize}
\end{assu}

Then the optimal control problem \eqref{promFBSC} (\autoref{prom-1}) can be reformulated as the following Stackelberg differential game problem.

\begin{prom}\label{prom-2}
  To find a map $\hat{\alpha}(\cdot):\mathcal{U}_{ad}[0,T]\rightarrow\mathcal{A}$ and a control $u^*(\cdot)$ such that
  \begin{equation}
      \begin{cases}
          J_1(u(\cdot),\hat{\alpha}(u(\cdot))(\cdot)) = \inf\limits_{(y_0,z(\cdot)) \in \mathcal{A}} J_1(u(\cdot),y_0,z(\cdot)) \qquad \forall u(\cdot)\in \mathcal{U}_{ad}[0,T], \vspace{1ex} \\
          J_2(u^*(\cdot),\hat{\alpha}(u^*(\cdot))(\cdot)) = \inf\limits_{u(\cdot)\in\mathcal{U}_{ad}[0,T]} J_2(u(\cdot),\hat{\alpha}(u(\cdot))(\cdot)).
      \end{cases}
  \end{equation}
  The state equation is given as \eqref{forFBSCsys}.
\end{prom}

Under \autoref{assu-3},
for any given $(u(\cdot), y_0, z(\cdot))$ and the initial state $a$,
there exists a strong solution $(x(\cdot),y(\cdot))$ to \eqref{forFBSCsys}.
If the above optimal pair $(\hat{\alpha}(\cdot),u^*(\cdot))$ in \autoref{prom-2} exists, we have $(y_0^*,z^*(\cdot)) = \hat{\alpha}(u^*(\cdot))(\cdot)$ and $u^*(\cdot)$ is an optimal control of \autoref{prom-1}. The advantage of \autoref{prom-2} is that it does not need the regularity/integrability on $z(\cdot)$,
and the state equation is a forward SDE instead of a backward SDE as in \autoref{prom-1}.
The difficulty for solving \autoref{prom-2} is that it has to treat another optimization goal $E \left[|y(T)-g(x(T))|^2\right]$.

For simplicity, we denote the leader's cost functional $J_2(u(\cdot),y_0,z(\cdot))$ by $J(u(\cdot))$,
and the follower's cost functional $J_1(u(\cdot),y_0,z(\cdot))$ by $J(y_0,z(\cdot))$ in the following of the paper.

\section{Deep neural network for solving the stochastic Stackelberg differential game problem}\label{sec-schemes}

In this section, we consider to solve \autoref{prom-2} with deep neural network.
As shown in \autoref{sec-form}, \autoref{prom-2} is an optimal control problem with forward state equations,
it contains two optimization goals,
and the follower's cost functional $J(y_0,z(\cdot))$ should be optimized in priority to the leader's cost functional.
In order to solve \autoref{prom-2},
instead of using the penalty method which multiply the follower's cost functional $J(y_0,z(\cdot))$ by a penalty parameter and adding it to the leader's cost functional $J(u(\cdot))$,
we propose a bi-level optimization method.
We optimize these two cost functionals alternatively and spend more computational cost on the follower's cost functional $J(y_0,z(\cdot))$.
In this way, on the one hand, \autoref{prom-2} become unconstrained with respect to the control.
On the other hand, we do not have to choose the approximate penalty parameter,
which is usually hard to choose in the penalty method.
The detailed process for optimizing these two functionals can be found in \autoref{ssec-UpPama}.

Before showing the approximation algorithm of \autoref{prom-2},
we firstly need to discretize the forward state equation \eqref{forFBSCsys}.
Let $N\in\mathbb{N}^+$ and $ 0 = t_{0}<t_{1}<t_{2}<\cdots<t_{N-1}<t_{N} =T $ so that
\begin{equation*}
    \delta = \max_{0\leq i \leq N-1} |t_{i+1}-t_i|,
\end{equation*}
is sufficiently small.
Define $ \Delta t_{i}=t_{i+1}-t_{i} $ and $ \Delta B_{t_i}=B(t_{i+1})-B(t_{i}) $,
where $ B(t_i) \sim \mathcal{N}(0,t_{i}) $, for $ i = 0, 1, 2,\cdots, N-1 $.
We apply the Euler-Maruyama scheme to the forward state equation \eqref{forFBSCsys},
then we have
\begin{equation}\label{dis-FBSC}
  \left\{
    \begin{array}{l}
      x_{t_{i+1}} = x_{t_i} + b(t_i,x_{t_i},y_{t_i},z_{t_i},u_{t_i}) \Delta t_i + \sigma(t_i,x_{t_i},y_{t_i},z_{t_i},u_{t_i}) \Delta B_{t_i}, \vspace{1ex} \\
      y_{t_{i+1}} = y_{t_i} - l(t_i,x_{t_i},y_{t_i},z_{t_i},u_{t_i}) \Delta t_i + z_{t_i} \Delta B_{t_i}, \vspace{1ex} \\
      x_{t_{0}} = a, \qquad y_{t_{0}} = y_0.
    \end{array}
    \right.
\end{equation}
For easier expression, we denote $x(t_i)$ as $x_{t_i}$ here.

We use Monte Carlo sampling to approximate the expectations in the cost functionals.
Then, the leader's cost functional \eqref{CostFB} can be evaluated by the following schemes:
\begin{equation}\label{eq-dis-obj-u}
    J(u(\cdot)) = \dfrac{1}{M} \sum_{m=1}^{M} \left[ \sum_{i=0}^{N-1} f(t_i,x_{t_i}^m,y_{t_i}^m,z_{t_i}^m,u_{t_i}^m) \Delta t_i + h(x_{t_N}^m) + \gamma(y_{t_0}^m) \right],
\end{equation}
where $M$ represents the number of Monte Carlo samples.
Similarly, the follower's cost functional $J(y_0,z(\cdot))$ can be approximated by
\begin{equation}\label{eq-dis-obj-z}
    J(y_0,z(\cdot)) = \dfrac{1}{M} \sum_{m=1}^{M} \left[ |y_{t_N}^m - g(x_{t_N}^m)|^2 \right],
\end{equation}
Here the numbers of Monte Carlo samples $M$ in \eqref{eq-dis-obj-u} and \eqref{eq-dis-obj-z} may be different.
Without causing confusion, we use the same notation $M$ to represent the Monte Carlo sampling number, and the same notation $J(u(\cdot))$ and $J(y_0,z(\cdot))$ as in the continuous form to represent the corresponding discrete cost functionals.

\subsection{Neural network architecture}\label{ssec-Netwok}

In this subsection, we present the neural network architecture for solving \autoref{prom-2} approximately.
Let $q\in\mathbb{N}^+$ be a sufficiently large natural number,
and $\mathcal{NN}^{\theta_u}:[0,T]\times\mathbb{R}^n \mapsto U$ be a Borel measurable function.

Now we let the function $\mathcal{NN}^{\theta_u}$ with suitable $\theta_u \in \mathbb{R}^q$ be the approximation of the control $u$:
\begin{align}
    u_{t_i} &= \mathcal{NN}^{\theta_u}(t_i, x_{t_i};\theta_u), \label{EqApp}
\end{align}
for $i\in\{0,1,\cdots,N-1\}$.
$\theta_u$ is the trainable parameter in this neural network.
We represent the function $\mathcal{NN}^{\theta_u}$ with a multilayer feed-forward neural network of the form
\begin{equation}\label{eqBaseNet}
    \mathcal{NN}^{\theta_u} = \psi \circ \mathcal{A}_{I} \circ \sigma_{I -1} \circ \mathcal{A}_{I -1} \circ \cdots \circ \sigma_{1} \circ \mathcal{A}_1,
\end{equation}
where
\begin{itemize}
    \item $I$ is a positive integer specifying the depth of the neural network,
    \item $\mathcal{A}_1,\cdots,\mathcal{A}_{I}$ are functions of the form
    \begin{equation*}
        \mathcal{A}_i = w_i \alpha_i + b_i \in \mathbb{R}^{d_{i}}, \qquad \text{for } 1\leq i\leq I, \qquad \alpha_i\in \mathbb{R}^{d_{i-1}},
    \end{equation*}
    where the matrix weights $w_i$ and the bias vector $b_i$ are trainable parameters such that $\theta_u=(w_i,b_i)_{1\leq i\leq I}$,
    $d_i$ is the number of nodes in layer $i$,
    and $\alpha_0$ represents the inputs of the neural network;
    \item $\sigma_{I-1},\cdots,\sigma_1$ are the nonlinear activation functions, such as sigmoid, ReLU, ELU, etc.;
    \item $\psi : \mathbb{R}^{d_I} \mapsto U$ is a given function, in this paper, we set $d_I=k$ and $U=\mathbb{R}^k$.
\end{itemize}

\begin{figure}[htbp]
  \centering
  \includegraphics[width=0.95\textwidth]{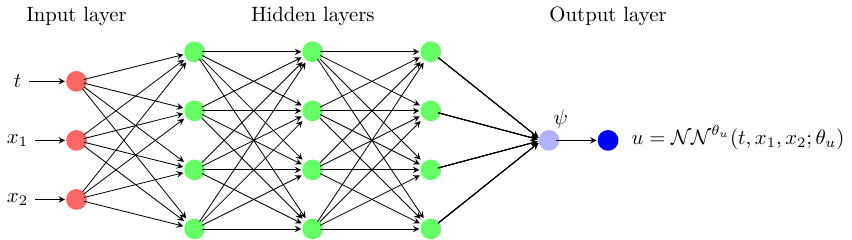}
  \caption{Representation of a single neural network with $I=4$, $d_1=d_2=d_3=4$, $n=2$, $d_4=k=1$.}
  \label{fig-BaseNetSingle}
\end{figure}

We treat the time $t$ as a part of input variables and the dimensions of the input and output in the neural network are $(n+1)$ and $k$, respectively.
\autoref{fig-BaseNetSingle} shows an example of a single neural network for $n=2$, $m=k=1$.

In this paper,
as the time $t$ is treated as an input variable,
we use a common neural network, i.e. we share the parameters of the neural network among all the time points.
The neural network is constructed with 5 layers,
including 1 input layer with $(n+1)$ neurons,
3 hidden layers with $(n+10)$ neurons and 1 output layer with $k$ neurons.
We adopt ReLU as the activation functions through the network, and add a batch normalization layer for each layer.

\begin{figure}[H]
  \centering
  \includegraphics[scale=1.1]{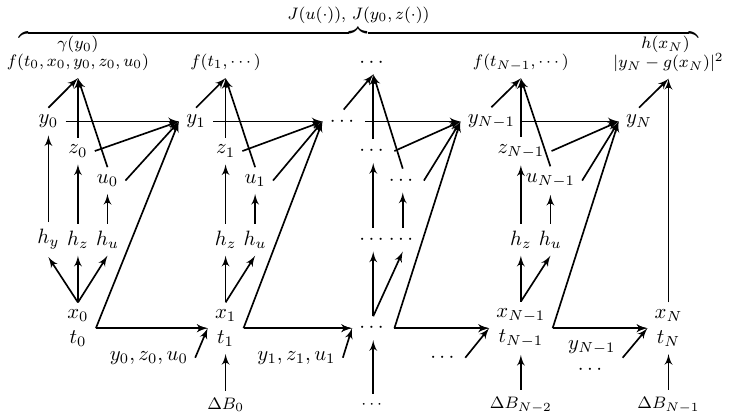}
  \caption{
  For easier notation, here we use $x_i$ to represent $x_{t_i}$, and $y_i$, $z_i$, $u_i$ for $y_{t_i}$, $z_{t_i}$, $u_{t_i}$ respectively.
  In the whole neural network structure,
  $h_y$, $h_z$ and $h_u$ represent the hidden layers of $\mathcal{NN}^{\theta_y}$, $\mathcal{NN}^{\theta_z}$ and $\mathcal{NN}^{\theta_u}$, respectively.
  And a common neural network is used for all the discrete time points.}
  \label{fig-NetConstruce}
\end{figure}

For approximating the controls $y_0$ and $z(\cdot)$,
we also construct two different feed-forward neural networks $\mathcal{NN}^{\theta_y}$ and $\mathcal{NN}^{\theta_z}$, respectively.
\begin{equation}\label{eq-simYZ}
  \left\{
  \begin{array}{l}
    y_0 = \mathcal{NN}^{\theta_y}(x_0;\theta_y) \vspace{1ex} \\
    z_{t_i} = \mathcal{NN}^{\theta_z}(t_i, x_{t_i};\theta_z).
  \end{array}
  \right.
\end{equation}
The architectures of $\mathcal{NN}^{\theta_y}$ and $\mathcal{NN}^{\theta_z}$ are similar to that of $\mathcal{NN}^{\theta_u}$,
and the trainable parameters of them are $\theta_y$ and $\theta_z$, respectively.
They have the same depth,
the same dimensions for the hidden layers, the same activation functions with the neural network $\mathcal{NN}^{\theta_u}$, but different input and output dimensions.
The dimensions of the input and output for $\mathcal{NN}^{\theta_y}$ are $n$ and $m$,
while that of the input and output for $\mathcal{NN}^{\theta_z}$ are $n+1$ and $m\times d$.
We also add a batch normalization layer for each layer in the two neural networks. The whole neural network structure is shown in \autoref{fig-NetConstruce}.

\subsection{Updating the parameters in the neural networks}\label{ssec-UpPama}

From the formulation of \autoref{prom-2},
two objective functionals $J(u(\cdot))$ and $J(y_0,z(\cdot))$ should be optimized, and the follower's cost functional $J(y_0,z(\cdot))$ should be optimized in priority to the leader's cost functional.
Which means that we should make much more effort on optimizing the follower's cost functional $J(y_0,z(\cdot))$.
A commonly used method for solving this problem is the penalty method,
the idea of which is to treat $J(y_0,z(\cdot))$ as a penalty and add it to the leader's cost functional $J(u(\cdot))$:
\begin{equation}\label{eqObjectAll}
    \bar{J}(y_0,z(\cdot),u(\cdot)) = J(u(\cdot)) + \mu J(y_0,z(\cdot)),
\end{equation}
where $\mu$ is a sufficiently large penalty parameter.
And the aim is to find the optimal control $(y^*_0,z^*(\cdot),u^*(\cdot))$ of \eqref{eqObjectAll},
in this way, \autoref{prom-2} is transformed to a classical stochastic optimal control problem without state constraint.

In this paper, different from the penalty method \eqref{eqObjectAll},
we develop a novel deep learning method for solving this kind of stochastic Stackelberg differential game problem.
In the neural network architecture, we update the network parameters through bi-level optimization of the two cost functionals \eqref{eq-dis-obj-u} and \eqref{eq-dis-obj-z}.
Specifically, for a given network parameter $\theta_u$,
we update the network parameters $(\theta_y, \theta_z)$ through the discrete follower's cost functional \eqref{eq-dis-obj-z}.
Then, under the approximation of $y_0,z(\cdot)$ with the updated $(\theta_y, \theta_z)$,
we update the network parameters $\theta_u$ through the discrete leader's cost functional \eqref{eq-dis-obj-u}. In other words, in each training step, the network parameters $\theta_u$ for optimizing \eqref{eq-dis-obj-u} and $(\theta_y, \theta_z)$ for optimizing \eqref{eq-dis-obj-z} are updated alternatively.
As mentioned above, the minimum of the follower's cost functional $J(y_0,z(\cdot))$ should be satisfied in priority to the leader's cost functional. Therefore in each training step, the update times of $(\theta_y, \theta_z)$ will be much higher than that of $\theta_u$ in practice.


Formally speaking, for a certain $\kappa\in\mathbb{N}^+$, a training step contains two different sub-steps
\begin{equation}\label{eqTheta}
    \left\{
    \begin{array}{l}
      \theta^{l+1}_u = \Xi_1 (\theta^l_u, \nabla J(u(\cdot))), \vspace{1ex} \\
      (\theta^{l+1}_y,\theta^{l+1}_z) = \Xi_2 \circ \cdots \circ \Xi_2 \circ \Xi_2 (\theta^{l}_y,\theta^{l}_z, \nabla J(y_0,z(\cdot))), \quad \mbox{total $\kappa$ times,}
    \end{array}
    \right.
\end{equation}
where $\Xi_1 (\theta^l_u, \nabla J(u(\cdot)))$ and $\Xi_2 (\theta^{l}_y,\theta^{l}_z, \nabla J(y_0,z(\cdot)))$ are given functions of parameters and their gradients.
They represent the parameter update operations for the cost functionals \eqref{eq-dis-obj-u} and \eqref{eq-dis-obj-z}, respectively.
When the parameter $\theta^{l}_u$ is updated through $\Xi_1 (\theta^l_u, \nabla J(u(\cdot)))$,
parameters $(\theta^{l}_y,\theta^{l}_z)$ remain unchanged,
and when parameters $(\theta^{l}_y,\theta^{l}_z)$ are updated through $\Xi_2^1 (\theta^{l}_y,\theta^{l}_z, \nabla J(y_0,z(\cdot)))$,
the parameter $\theta^{l}_u$ remains unchanged.
Here we call $\kappa$ the penalty updating coefficient,
which represents the ratio of the number of optimizations of $(\theta^{l}_y,\theta^{l}_z)$ to that of $\theta^{l}_u$.

In the actual calculation process,
since the updates of $\theta^l_u$ and $(\theta^{l}_y,\theta^{l}_z)$ are carried out alternatively,
it is feasible to perform the first update for $\theta^l_u$ or $(\theta^{l}_y,\theta^{l}_z)$.
And we only need to ensure that $\kappa$ is large.
The pseudo code of the proposed algorithm is given in Algorithm \ref{alg-main}.

\begin{algorithm}[H]
  \renewcommand{\thealgorithm}{1}
  \caption{The bi-level optimization algorithm for solving stochastic control problem driven by fully coupled FBSDE}
  \begin{algorithmic}[1]
    \Require The Brownian motion $ \Delta B_{t_i} $, initial state $a$ and parameter $ (\theta_y^0,\theta_z^0,\theta_u^0) $ and hyper-parameter $\kappa'$, etc.;
    \Ensure The precesses $ (x^{\ell}_{t_i},y^{\ell}_{t_i},z^{\ell}_{t_i},u^{\ell}_{t_i}) $, etc.
    \For { $ \ell = 0 $ to $ maxstep $}
    \State $ x_{0}^{\ell,m} = a $, $loss = 0$;
    \State $ y_{0}^{\ell,m} = \mathcal{NN}^{\theta_y}(x_{0}^{\ell,m};\theta_y^{\ell}) $;
    \State $loss = loss + \gamma(y_{0}^{\ell,m})$
    \For { $ i = 0 $ to $ N-1 $}
    \State $u^{\ell,m}_{t_i} = \mathcal{NN}^{\theta_u}(t_i, x^{\ell,m}_{t_i};\theta_u^{\ell});$
    \State $z^{\ell,m}_{t_i} = \mathcal{NN}^{\theta_z}(t_i, x^{\ell,m}_{t_i};\theta_z^{\ell});$
    \State $x^{\ell,m}_{t_{i+1}} = x^{\ell,m}_{t_i} + b(t_i,x^{\ell,m}_{t_i},y^{\ell,m}_{t_i},z^{\ell,m}_{t_i},u^{\ell,m}_{t_i})\Delta t_i + \sigma(t_i,x^{\ell,m}_{t_i},y^{\ell,m}_{t_i},z^{\ell,m}_{t_i},u^{\ell,m}_{t_i}) \Delta B_{t_i};$
    \State $y^{\ell,m}_{t_{i+1}} = y^{\ell,m}_{t_i} - l(t_i,x^{\ell,m}_{t_i},y^{\ell,m}_{t_i},z^{\ell,m}_{t_i},u^{\ell,m}_{t_i})\Delta t_i + z^{\ell,m}_{t_i} \Delta B_{t_i};$
    \State $loss = loss + f(t_i,x^{\ell,m}_{t_i},y^{\ell,m}_{t_i},z^{\ell,m}_{t_i},u^{\ell,m}_{t_i})\Delta t_i;$
    \EndFor
    \State $J(u(\cdot)) = \dfrac{1}{M} \sum\limits_{m=1}^M \left[loss + h(x_{t_N}^{\ell,m})\right];$
    \State $J(y_0,z(\cdot)) = \dfrac{1}{M} \sum\limits_{m=1}^M \left[|y_{t_N}^{\ell,m}-g(x_{t_N}^{\ell,m})|^2\right];$
    \If {$\ell \mbox{ mod } (\kappa+1)=0$}
    \State $ \theta^{\ell}_u = \Xi_1 (\theta^{\ell}_u, \nabla J(u(\cdot))); $
    \Else
    \State $ (\theta^{\ell}_y,\theta^{\ell}_z) = \Xi_2 (\theta^{\ell}_y,\theta^{\ell}_z, \nabla J(y_0,z(\cdot))); $
    \EndIf
    \State $ (\theta_y^{\ell+1},\theta_z^{\ell+1},\theta_u^{\ell+1}) = (\theta_y^{\ell},\theta_z^{\ell},\theta_u^{\ell}); $
    \EndFor
    \State $J(u(\cdot)) = \dfrac{1}{M} \sum\limits_{m=1}^M \left[\sum\limits_{i=0}^{N-1} f(t_i,x^{\ell,m}_{t_i},y^{\ell,m}_{t_i},z^{\ell,m}_{t_i},u^{\ell,m}_{t_i})\Delta t_i + h(x_{t_N}^{\ell,m}) + \gamma(y_{0}^{\ell,m})\right].$
  \end{algorithmic}
  \label{alg-main}
\end{algorithm}

In Algorithm \ref{alg-main},
$maxstep$ represents the total number of parameter update times. That is to say, we perform $maxstep/(\kappa+1)$ training steps. And in each training step,
there are $\kappa$ times updates for the parameters $(\theta_y,\theta_z)$ and one time update for the parameter $\theta_u$.
In this bi-level optimization method, the choices of $\Xi_1$ and $\Xi_2$ are independent,
and the training parameters, such as the learning rate, the optimizer, the neural network architecture, can be chosen separately.
We can choose an appropriate penalty updating coefficient $\kappa$ in $\Xi_2$ to make the follower's cost functional $\eqref{eq-dis-obj-z}$ small enough,
instead of simply choosing the approximate penalty coefficient $\mu$ in the penalty method.

The above proposed algorithm can be regarded as a relaxing method for solving \autoref{prom-2}.
In \autoref{prom-2}, the optimization of the leader's cost functional $J(u(\cdot))$ should be based on that of the follower's cost functional $J(y_0,z(\cdot))$, thus we should firstly guarantee that $ E \left[|y(T)-g(x(T))|^2\right]=0 $. In our proposed method, we relax this constraint to improve the computation efficiency. We make the value of $ \left[E |y(T)-g(x(T))|^2\right] $ sufficiently close to 0, but not strictly equal to 0. Specifically, we make the coefficient $\kappa$ to be large enough so that we can obtain an enough small value of the follower's cost functional.

\section{Numerical results}\label{sec-numerical-results}

In this section,
we show two optimal investment-consumption portfolio examples solved through recursive utility with our proposed algorithm.
If not specially mentioned, we use a 5-layer fully connected neural network, 512 samples of Brownian motion in the test set, and the number of time points is $N=25$.
The implementations are performed through \textsc{TensorFlow} on a \textsc{Lenovo} computer with a 2.40 Gigahertz (GHz) Inter Core i7 processor and 8 gigabytes (GB) random-access memory (RAM).

Here we introduce a continuous and stochastic recursive utility problem which was first introduced by Duffie and Epstein \cite{duffie1992stochastic}.
Suppose there are $n+1$ assets trading continuously in the market,
one of which is a risk-free asset whose price process is given as
\begin{equation}\label{eqRiskFreeAsset}
    d P^0(t) = r(t)P^0(t) dt, \qquad P(0)^0 = p^0 >0,
\end{equation}
where $r(t)$ is the instantaneous rate of return.
The other $n$ assets are risky assets satisfying
\begin{equation}\label{eqRiskyAssets}
    d P^i(t) = P^i(t) \left\{ \mu^i(t) dt + \left\langle \sigma^i(t), dB(t) \right\rangle \right\}, \qquad P(0)^i = p_0^i>0,
\end{equation}
where $\mu^i(t):[0,T]\times\Omega\mapsto\mathbb{R}$ is the appreciation rate,
$\sigma^i(t):[0,T]\times\Omega\mapsto\mathbb{R}^n$ is the volatility of the stocks and $B(t)$ is a standard Brownian motion valued in $\mathbb{R}^n$.
All the processes are assumed to be $\mathcal{F}_t$-adapted.

An investor starts with a given initial wealth $x_0$ and his total wealth at time $t$ is denoted as $x(t)$.
The wealth process $\{x(t)\}_{0\leq t\leq T}$ is given by
\begin{equation}\label{eqWealthPro}
    \left\{
    \begin{array}{l}
        dx(t) = \left\{ x(t) r(t) + \sum\limits_{i=1}^n x(t)\pi^i(t)(\mu^i(t)-r(t)) - c(t) \right\} dt + \sum\limits_{i=1}^n x(t)\pi^i(t) \left\langle \sigma^i(t), dB(t) \right\rangle, \vspace{1ex} \\
        x(0)=x_0>0,
    \end{array}
    \right.
\end{equation}
where $\pi(t) = (\pi^1(t), \pi^2(t), \cdots, \pi^n(t))$ is called the portfolio of the investor and $c(t)$ is called the consumption plan process.
$\pi^i(t)$ represents the proportion of total wealth $x(t)$ invested in the $i$-th risky asset and is taking value in $[0,1]$,
which means that the short-selling is prohibited.
The remaining fraction $\pi^0(t) = 1- \sum_{i=1}^n \pi^i(t) $ valued in $[0,1]$ is thus the proportion of the left in the form of risk-free bond.
In addition, we assume that the consumption process is non-negative ($c(t) \geq 0)$.

The utility of the investor at time $t$ is denoted as a function of the instantaneous consumption $c(t)$ and the future utility.
The utility process can be regarded as the solution of a BSDE given as
\begin{equation}\label{eqRecurBSDE}
    -dy(t) = l(t,c(t),y(t),z(t))dt - \left\langle z(t), dB(t) \right\rangle, \qquad y(T) = g(x(T)),
\end{equation}
where $g(x)$ is the terminal utility function.
Under some regularity assumptions (such as \autoref{assu-1}, \autoref{assu-2}),
equation \eqref{eqRecurBSDE} has a unique solution $(y(\cdot),z(\cdot))$.
More precisely, the recursive utility at time $t$ can be denoted by
\begin{equation}\label{eqRecurUtility}
    y(t) = E \left[\int_t^T l(s,c(s),y(s),z(s)) ds + g(x(T)) \Big | \mathcal{F}_t\right] ,
\end{equation}
where $\mathcal{F}_t$ represents the natural filtration associated with the standard Brownian motion.

The goal of the investor is to choose the optimal investment and consumption to maximize the utility at time zero.
Mathematically, he aims to find the optimal controls $\pi^*(\cdot)$ and $c^*(\cdot)$,
such that
\begin{equation}\label{eqRecurProm}
    y^*(0) = \sup_{(\pi(\cdot),c(\cdot))\in\mathcal{U}_{ad}[0,T]} y_0,
\end{equation}
where $\mathcal{U}_{ad}[0,T]$ represents the set of all possible investment-consumption strategies and $y_0$ is the recursive utility given as \eqref{eqRecurUtility}.
We rewrite $y_0$ as
\begin{equation}\label{eqY0Target}
    y_0 = J(\pi(\cdot),c(\cdot)) = E \left[\int_0^T l(t,c(t),y(t),z(t)) dt + g(x(T))\right] ,
\end{equation}
and we suppose that $l$ is a generalized recursive utility function given as:
\begin{equation}\label{eqGenRecDyn}
    \begin{aligned}
        l(t,c,y,z) = -\beta y - w(z) + u(c).
    \end{aligned}
\end{equation}
Here $u(c)$ is the instantaneous utility function.
The function $w(z)$ is a differentiable function and can be regarded as a risk-aversion coefficient.

Actually, the recursive utility problem \eqref{eqRecurProm} is a special case of Problem \ref{prom-1},
in which $f(t,x,y,z,u) = 0$, $h(x)=0$ and $\gamma(y)=y$,
therefore we can reformulate it to Problem \ref{prom-2}, where the goal of the leader is to find the optimal controls $(\pi^*(t),c^*(t))_{0\leq t\leq T}$ which maximize \eqref{eqY0Target}, and the goal of the follower is to minimize
\begin{equation}\label{eqRecOptFBSDE}
    J(y_0,z(\cdot)) = E \left[|y(T) - g(x(T))|^2\right] .
\end{equation}
under the given controls $(\pi^*(t),c^*(t))$. Then we solve the problem  through our bi-level optimization algorithm presented in \autoref{sec-schemes}.

In this example, in order to approximate the controls $\pi(t)$ and $c(t)$, we need to construct two neural networks $\mathcal{NN}_{\pi}$, $\mathcal{NN}_c$ , and both of the controls are supposed as the feedback controls of the time and the wealth process $(t,x(t))$.
The neural network approximating $\pi(t)$ contains one $(n+1)$-dimensional input layer,
three  $100$-dimensional hidden layers and one $(n+1)$-dimensional output layer.
For the constraint of the portfolio,
$\pi^i(t)\in[0,1]$ and $\sum_{i=0}^n \pi^i(t) =1$ for $i=0,\cdots,n$,
we deal with it through a softmax function $\psi$.
The function $\psi$ is given as
\begin{equation*}
    \psi(x) = \left(\dfrac{e^{x^0}}{\sum_{i=0}^n{e^{x^i}}}, \dfrac{e^{x^1}}{\sum_{i=0}^n{e^{x^i}}}, \cdots, \dfrac{e^{x^n}}{\sum_{i=0}^n{e^{x^i}}}\right), \qquad x=(x^0,\cdots,x^n) \in \mathbb{R}^{n+1},
\end{equation*}
which is often used in classification problems such as image classification. Intuitively, it makes sense to choose the softmax function because choosing the best investment asset based on the current state is essentially a classification problem. Moreover, we use a function $\psi(x) = \max(x, 0)$ to deal with the constraint $c(t)\geq 0$. In addition, we construct the other two neural networks to approximate $y_0$ and $z(t)$, respectively. In all, four feed-forward neural networks are constructed in this example, the first for approximating the consumption rate $c(t)$, the second for approximating the wealth proportion $\pi(t)$, the third for approximating $y_0$, and the last for approximating $z(t)$.

Note that in \eqref{eqRecOptFBSDE},
$y_0$ is regarded as a control parameter which is approximated through deep neural network. While at the same time, it is the optimization objective in \eqref{eqY0Target}. In order to make difference,
we call the $y_0$ in \eqref{eqY0Target} \textbf{the integral form} and that in the neural network \textbf{the parametric form}.
And we should emphasize that instead of optimizing the parametric form $y_0$, we optimize the integral form $y_0 = J(\pi(\cdot),c(\cdot))$ with respect to $(\pi(\cdot),c(\cdot))$. The reason is that the derivatives of the parametric form $y_0$ with respect to $(\pi(\cdot),c(\cdot))$ do not exist. What's more, we compare the integral form with the parametric form of $y_0$. We measure the distance of the values between the parametric form  and the integral form and take it as a criterion for the effectiveness of our proposed algorithm.

\subsection{Case 1: a linear driver }

In the first example,
we consider a linear driver and set $w(z) \equiv 0$,
then the recursive utility functional we need to maximize is given as
\begin{equation}\label{eqY0Target-y}
    y_0 = J(\pi(\cdot),c(\cdot)) = E \left[\int_0^T (-\beta y(t) + u(c(t))) dt + g(x(T))\right] .
\end{equation}
In this case, \eqref{eqY0Target-y} is equivalent  to the stochastic functional
\begin{equation}\label{eqY0Target-1}
    y_0 = J'(\pi(\cdot),c(\cdot)) = E \left[\int_0^T e^{-\beta t}u(c(t)) dt + e^{-\beta T} g(x(T))\right] .
\end{equation}
Thus the stochastic recursive utility problem \eqref{eqWealthPro} and \eqref{eqY0Target} can degenerate to a classic stochastic optimal control problem as \eqref{eqWealthPro} and \eqref{eqY0Target-1}.
There have been some deep learning methods for solving the stochastic control problem \eqref{eqWealthPro}-\eqref{eqY0Target-1},
such as the methods in \cite{han2016deep,pham2018deep_2,ji2021solving}.

Here we suppose that $u(c)$ is a quadratic utility function
\begin{equation*}
    u(c) = (c-c^2),
\end{equation*}
and the terminal utility function $g(x)$ is
\begin{equation*}
    g(x) = e^{-x}.
\end{equation*}

We calculate the recursive utility problem with the previously mentioned investment-consumption model.
Let the wealth process and trading assets processes satisfy \eqref{eqRiskFreeAsset}-\eqref{eqWealthPro}. We set $r(t)=0.03$, $\mu^i(t) = 0.05$, $\sigma=(\sigma^i,\cdots,\sigma^n)=0.1I_n$ for $i=1,\cdots, n$,
and the initial wealth $x_0$ is equal to $100$.
The parameters $\beta$ in $l(t,c,y,z)$ is set to be $0.05$.
As a comparison,
we also calculate the classic stochastic optimal control problem \eqref{eqWealthPro} and \eqref{eqY0Target-1}.
Here we construct a similar neural network with that in \cite{han2016deep} and regard the discrete time point $t$ as an input of the neural network.

The implementation results with different dimensions and different terminal time are shown in \autoref{tab-Utility-y}.
All the results are take from 5 independent runs. The number of iterations for parameter update is 18,000 and $\kappa$ is set to be 19. That is to say, after 19 iterations of optimization for the follower's cost functional \eqref{eqRecOptFBSDE},
we perform one iteration of optimization for the leader's cost functional \eqref{eqY0Target}. And the total number of training steps is $18,000/(19+1)=900$.
We show the integral form and the parametric form of $y_0$ (Inte. $y_0$ and Para. $y_0$ in the table, respectively),
and calculate their mean values and variances among 5 independent runs.
We also show the values calculated through the classic stochastic control problem \eqref{eqWealthPro} and \eqref{eqY0Target-1} (Clas. $y_0$ in the table).

\begin{table}[htbp]
  \centering
  \captionsetup{font={small}}
  \caption{The implementation results of $y_0$ for case 1}
  \resizebox{\textwidth}{!}{
  \begin{tabular}{|c|c|c|c|c|c|c|c|c|c|}
    \hline
        \multicolumn{2}{|c|}{} & \multicolumn{2}{|c|}{$T=0.25$}&  \multicolumn{2}{|c|}{$T=0.50$}&  \multicolumn{2}{|c|}{$T=0.75$}&  \multicolumn{2}{|c|}{$T=1.00$} \\
    \hline
    \multicolumn{2}{|c|}{}& Mean& Var.& Mean&Var.& Mean&Var.& Mean&Var.\\
    \hline
    \multirow{3}{*}{$n=10$} & Clas. $y_0$ & 0.06192 & 8.0998e-08 & 0.12323 & 7.2303e-08 & 0.18416 & 1.7761e-10 & 0.24314 & 3.9102e-07 \\
    \cline{2-10}
        & Inte. $y_0$                     & 0.06161 & 2.6347e-08 & 0.12269 & 9.8344e-08 & 0.18299 & 2.9634e-07 & 0.24216 & 7.3878e-07 \\
    \cline{2-10}
        & Para. $y_0$                     & 0.06146 & 3.1114e-07 & 0.12194 & 3.0065e-07 & 0.18217 & 3.3235e-06 & 0.24137 & 1.4374e-06 \\
    \hline
    \multirow{3}{*}{$n=20$} & Clas. $y_0$ & 0.06183 & 1.2111e-07 & 0.12301 & 3.9042e-07 & 0.18389 & 1.1229e-07 & 0.24359 & 3.3391e-07 \\
    \cline{2-10}
        & Inte. $y_0$                     & 0.06169 & 2.3650e-08 & 0.12268 & 2.9813e-07 & 0.18248 & 3.5300e-07 & 0.24202 & 4.7493e-07 \\
    \cline{2-10}
        & Para. $y_0$                     & 0.06169 & 2.9493e-07 & 0.12238 & 8.1997e-07 & 0.18149 & 1.0537e-06 & 0.24084 & 5.2556e-06 \\
    \hline
    \multirow{3}{*}{$n=30$} & Clas. $y_0$ & 0.06160 & 3.7891e-07 & 0.12282 & 6.4121e-07 & 0.18376 & 7.0530e-08 & 0.24315 & 8.0416e-07 \\
    \cline{2-10}
        & Inte. $y_0$                     & 0.06137 & 4.5255e-08 & 0.12203 & 1.9316e-07 & 0.18204 & 1.0402e-07 & 0.24268 & 6.7658e-07 \\
    \cline{2-10}
        & Para. $y_0$                     & 0.06036 & 3.8921e-07 & 0.12169 & 5.4307e-06 & 0.17921 & 4.2815e-07 & 0.24098 & 5.2413e-06 \\
    \hline
    \multirow{3}{*}{$n=40$} & Clas. $y_0$ & 0.06171 & 1.1584e-07 & 0.11055 & 1.3581e-03 & 0.18333 & 4.5646e-07 & 0.24321 & 8.7775e-07 \\
    \cline{2-10}
        & Inte. $y_0$                     & 0.06142 & 2.0930e-08 & 0.12226 & 2.3594e-07 & 0.18235 & 1.6639e-07 & 0.24228 & 5.0222e-07 \\
    \cline{2-10}
        & Para. $y_0$                     & 0.06129 & 4.6174e-07 & 0.12092 & 2.7358e-06 & 0.18113 & 3.5940e-06 & 0.23994 & 4.3413e-06 \\
    \hline
    \multirow{3}{*}{$n=50$} & Clas. $y_0$ & 0.06170 & 9.4358e-08 & 0.12303 & 2.7746e-07 & 0.18324 & 1.2237e-07 & 0.24126 & 1.9552e-06 \\
    \cline{2-10}
        & Inte. $y_0$                     & 0.06143 & 2.0046e-08 & 0.12222 & 6.9066e-08 & 0.18273 & 5.7814e-08 & 0.24134 & 3.6883e-07 \\
    \cline{2-10}
        & Para. $y_0$                     & 0.06077 & 2.7034e-06 & 0.12160 & 3.2163e-06 & 0.18115 & 1.3218e-06 & 0.23722 & 1.3412e-06 \\
    \hline
  \end{tabular}
  }
  \label{tab-Utility-y}
\end{table}

\begin{figure}[htbp]
    \centering
    \includegraphics[scale=0.45]{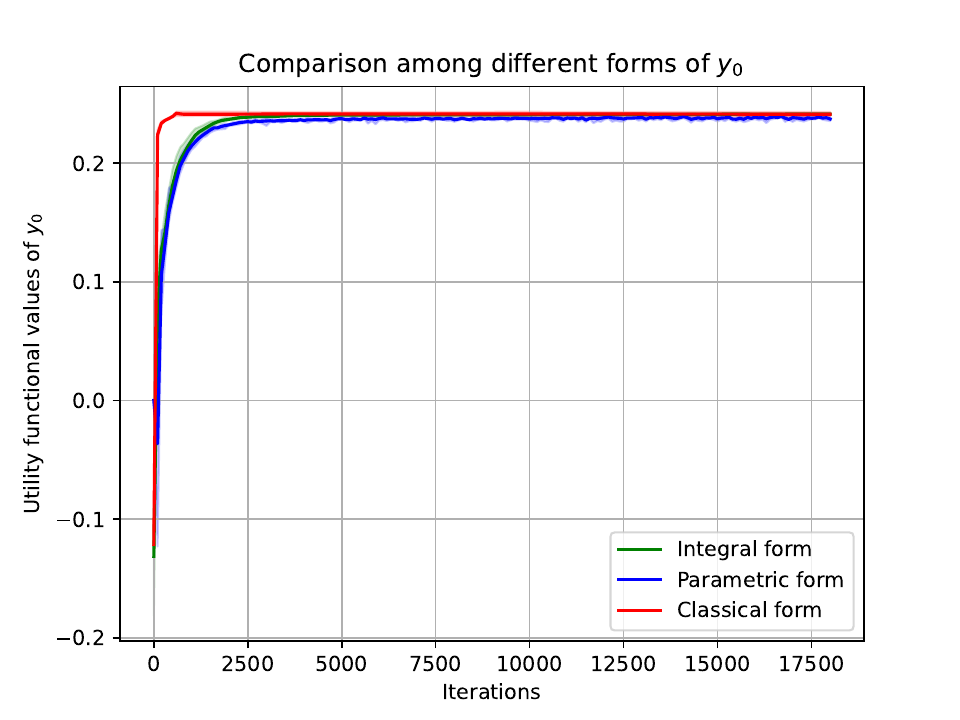}
    \includegraphics[scale=0.45]{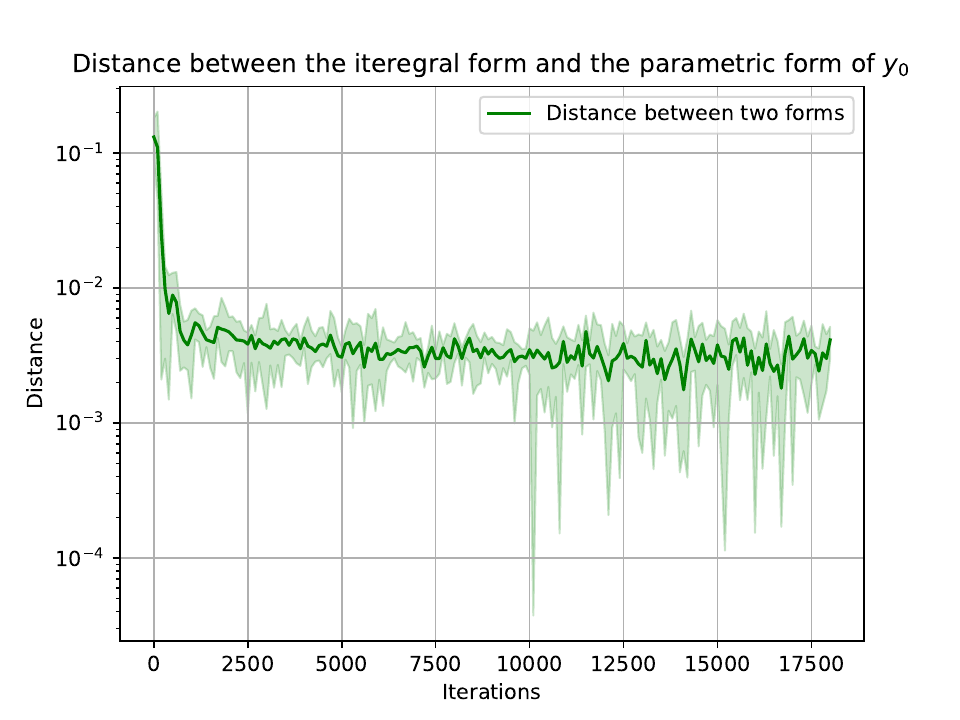}
    \caption{Case $n=50$ and $T=1.0$}
    \label{fig-Utility-y}
\end{figure}

The results in Table~\ref{tab-Utility-y} show that our proposed method is effective for the recursive utility problem.
The values of the utility $y_0$ between the integral form and the parametric form are close, and both of the two forms are close to the value of the classic form.
We also exhibit the curves of the utility values of different forms in Figure~\ref{fig-Utility-y}.
The left figure shows the mean and scope for the different forms of $y_0$ among 5 independent runs. We can see that the values of the three different forms of $y_0$ are getting closer with the increase of the number of iterations. At the same time, the value of the integral form is getting larger which meets our objective to maximize the recursive utility functional \eqref{eqY0Target} (see the green curve and scope). The right figure shows the mean and scope of the distance between the integral form and the parametric form.

\begin{figure}[htbp]
    \centering
    \includegraphics[scale=0.45]{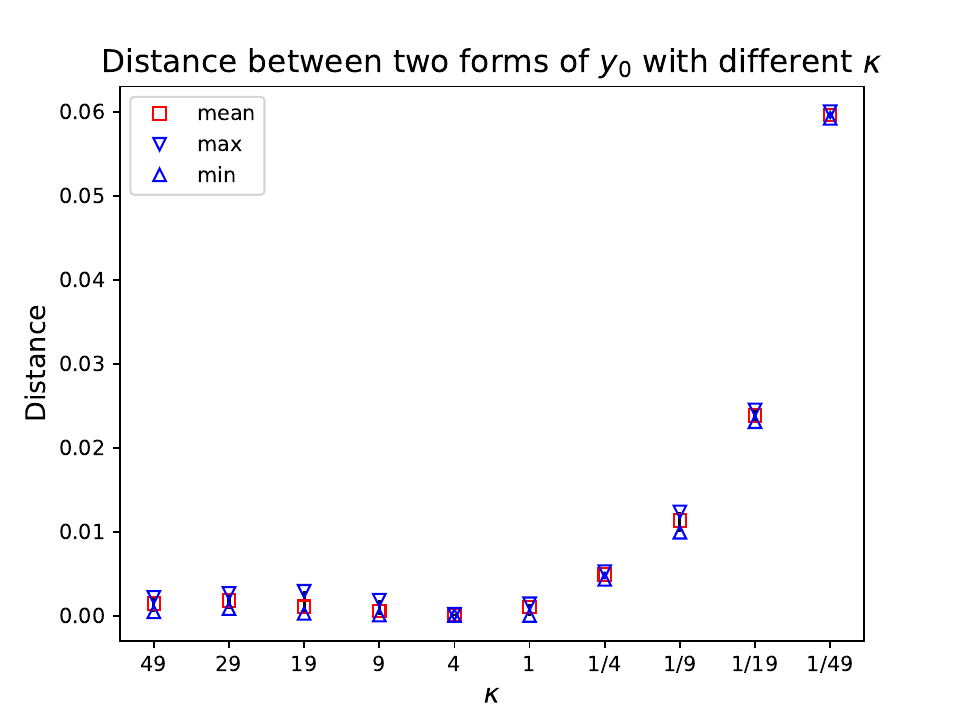}
    \includegraphics[scale=0.45]{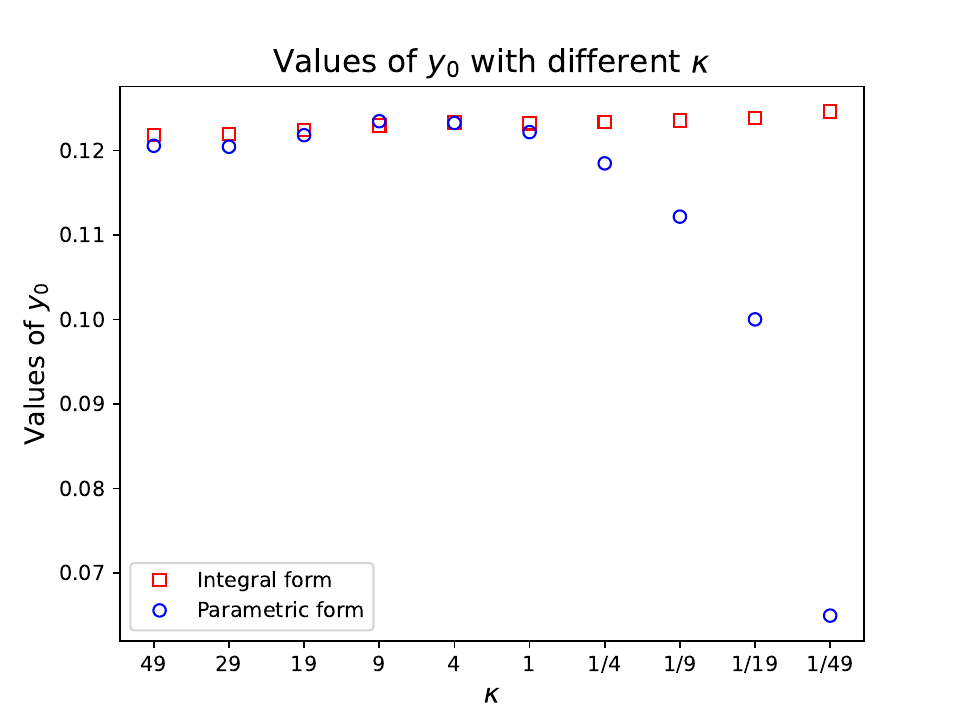}
    \caption{The distance of $y_0$ between the integral form and the parametric form and the values of $y_0$ with different $\kappa$ for $n=10$ and $T=0.5$.}
    \label{fig-UtilityModel-k-y}
\end{figure}

In order to show the impact of the penalty updating coefficient $\kappa$ to the convergence of the algorithm,
we vary the value of $\kappa$ for $T=0.5$ and $n=10$.
The results with different $\kappa$ are shown in \autoref{tab-Utility-k-y} and \autoref{fig-UtilityModel-k-y}.
In Figure~\ref{fig-UtilityModel-k-y}, the left figure shows that the distance between the integral form and the parametric form are small when $\kappa > 1$,
which means that the value of \eqref{eqRecOptFBSDE} can be small enough when we perform more iterations on \eqref{eqRecOptFBSDE}.
The right figure shows the  values of the utility functional are more stable for the case $\kappa > 1$ than the case $\kappa < 1$.
Here $\kappa < 1$ means that we perform more iterations of optimization on \eqref{eqY0Target} than that on \eqref{eqRecOptFBSDE}.

\begin{table}[htbp]
  \centering
  \caption{The implementation results of $y_0$ with different $\kappa$ for case 1}
  \begin{tabular}{|c|c|c|c|c|c|c|}
    \hline
    \multicolumn{1}{|c|}{}& \multicolumn{2}{|c|}{Inte. $y_0$} & \multicolumn{2}{|c|}{Para. $y_0$} & \multicolumn{2}{|c|}{Distance} \\
    \hline
    & Mean & Var. & Mean & Var. & Mean & Var. \\
    \hline
    $\kappa=49$   & 0.12184 & 1.9438e-07 & 0.12056 & 1.9288e-06 & 0.00145 & 5.2404e-07 \\
    \hline
    $\kappa=29$   & 0.12194 & 1.4726e-07 & 0.12045 & 2.4452e-06 & 0.00182 & 4.1077e-07 \\
    \hline
    $\kappa=19$   & 0.12243 & 7.4394e-08 & 0.12182 & 6.7248e-07 & 0.00107 & 2.2293e-07 \\
    \hline
    $\kappa=9$    & 0.12295 & 4.0809e-07 & 0.12347 & 1.2084e-08 & 0.00055 & 4.4437e-07 \\
    \hline
    $\kappa=4$    & 0.12332 & 6.8988e-09 & 0.12326 & 2.2820e-09 & 0.00010 & 5.4356e-09 \\
    \hline
    $\kappa=1$    & 0.12323 & 9.5936e-08 & 0.12218 & 5.6867e-08 & 0.00106 & 2.7134e-07 \\
    \hline
    $\kappa=1/4$  & 0.12340 & 1.9741e-08 & 0.11848 & 1.3041e-07 & 0.00492 & 1.3164e-07 \\
    \hline
    $\kappa=1/9$  & 0.12354 & 3.0263e-08 & 0.11218 & 6.6298e-07 & 0.01136 & 7.8541e-07 \\
    \hline
    $\kappa=1/19$ & 0.12386 & 2.7079e-08 & 0.10002 & 3.1441e-07 & 0.02384 & 2.6049e-07 \\
    \hline
    $\kappa=1/49$ & 0.12461 & 2.8228e-08 & 0.06496 & 9.7174e-08 & 0.05965 & 1.1475e-07 \\
    \hline
  \end{tabular}
  \label{tab-Utility-k-y}
\end{table}

\subsection{Case 2: a nonlinear driver}

In this case, the function $w(z)$ is given as $w(z) = \frac{\nu}{2} |z|^2$, $|z|^2=\left\langle z,z \right\rangle$, for $z\in \mathbb{R}^n$,
where $\nu$ is a given constant $\nu = 10$.
The functions $u(c)$, $g(x)$ and the other settings are the same as the linear driver.

\begin{figure}[htbp]
    \centering
    \includegraphics[scale=0.45]{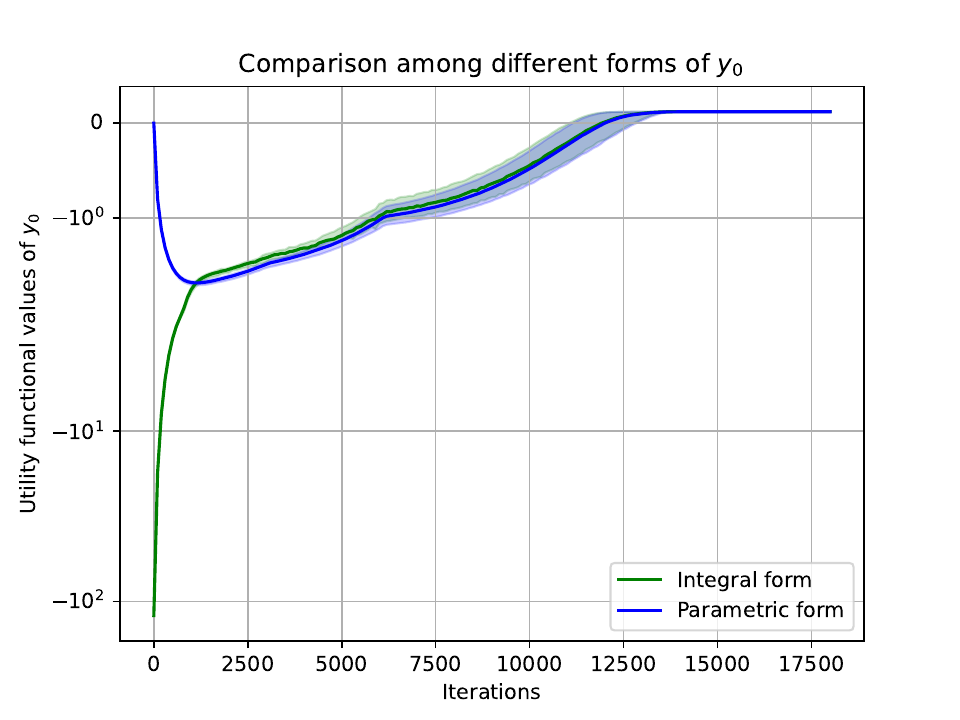}
    \includegraphics[scale=0.45]{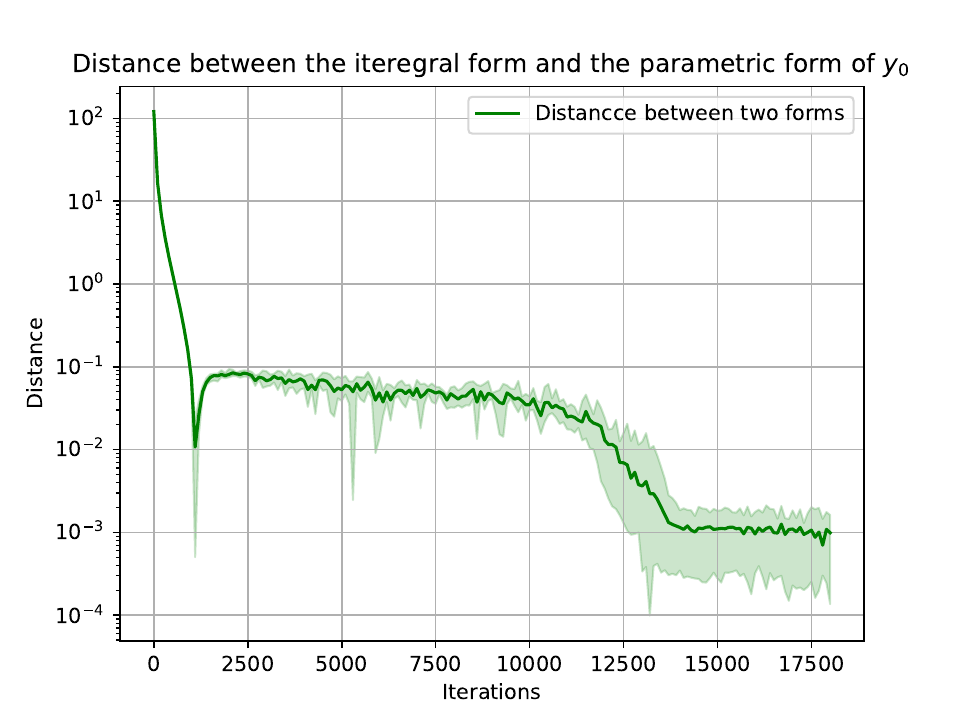}
    \caption{Case $n=50$ and $T=0.5$. The left figure shows the mean and scope for the values of $y_0$ among 5 independent runs. We can see that the two different forms of $y_0$ are getting closer when the number of iteration steps increases. Moreover, the value of the integral form of $y_0$ (see the green curve and scope) is getting larger and then tends to be stable. The right figure shows the mean and scope of the distances between the different forms of $y_0$.}
    \label{fig:UtilityModel}
\end{figure}

The implementation results with different dimensions and different terminal times for this example are shown in \autoref{tab:UtilityModel}.
We set $\kappa=19$  and the total number of iterations to be 18,000.
The computation results are similar to that with the linear driver. The distances between the integral form and the parametric form of $y_0$ are small at the end of the training. In \autoref{fig:UtilityModel}, we show the curves for the values and the distances of $y_0$ between the integral form and the parametric form with $n=50$ and $T=0.5$. We can see that the distances between the integral form and the parametric form of $y_0$ are getting closer to 0 when the number of iterations increases. Besides, the value of the integral form is getting larger, which shows that we can find the optimal investment-consumption strategy $(\pi^*(\cdot),c^*(\cdot))$ to maximize the recursive utility functional \eqref{eqY0Target}.

\begin{table}[htbp]
  \centering
  \captionsetup{font={small}}
  \caption{The implementation results of $y_0$ for case 2}
  \resizebox{\textwidth}{!}{
  \begin{tabular}{|c|c|c|c|c|c|c|c|c|c|}
    \hline
        \multicolumn{2}{|c|}{} & \multicolumn{2}{|c|}{$T=0.25$}&  \multicolumn{2}{|c|}{$T=0.50$}&  \multicolumn{2}{|c|}{$T=0.75$}&  \multicolumn{2}{|c|}{$T=1.00$} \\
    \hline
    \multicolumn{2}{|c|}{}& Mean& Var.& Mean&Var.& Mean&Var.& Mean&Var.\\
    \hline
    \multirow{3}{*}{$n=10$} & Inte. $y_0$ & 0.06166 & 2.2275e-08 & 0.12293 & 7.4394e-08 & 0.18284 & 4.2424e-07 & 0.24182 & 1.0601e-06 \\
    \cline{2-10}
        & Para. $y_0$                     & 0.06110 & 4.1587e-07 & 0.12240 & 6.7248e-07 & 0.18236 & 6.0408e-07 & 0.23978 & 8.4018e-06 \\
    \cline{2-10}
        & Distance                      & 0.00074 & 4.6313e-08 & 0.00064 & 2.2293e-07 & 0.00063 & 1.3188e-07 & 0.00236 & 2.5190e-06 \\
    \hline
    \multirow{3}{*}{$n=20$} & Inte. $y_0$ & 0.06156 & 1.3603e-07 & 0.12266 & 1.6195e-07 & 0.18284 & 9.5716e-08 & 0.24167 & 2.7468e-06 \\
    \cline{2-10}
        & Para. $y_0$                     & 0.06123 & 1.2727e-06 & 0.12302 & 3.9956e-07 & 0.18352 & 8.2581e-07 & 0.24052 & 1.1217e-05 \\
    \cline{2-10}
        & Distance                      & 0.00066 & 3.2714e-07 & 0.00054 & 7.1423e-08 & 0.00117 & 9.8714e-08 & 0.00164 & 1.7542e-06 \\
    \hline
    \multirow{3}{*}{$n=30$} & Inte. $y_0$ & 0.06157 & 4.3009e-08 & 0.12270 & 8.5964e-08 & 0.18256 & 8.1253e-07 & 0.24257 & 7.8675e-07 \\
    \cline{2-10}
        & Para. $y_0$                     & 0.06126 & 4.5989e-07 & 0.12206 & 1.2343e-06 & 0.18123 & 4.9427e-06 & 0.24145 & 2.5426e-06 \\
    \cline{2-10}
        & Distance                      & 0.00052 & 5.2667e-08 & 0.00109 & 3.3751e-07 & 0.00155 & 1.1978e-06 & 0.00112 & 5.7840e-07 \\
    \hline
    \multirow{3}{*}{$n=40$} & Inte. $y_0$ & 0.06145 & 2.8057e-08 & 0.12221 & 3.3040e-07 & 0.18304 & 6.4776e-07 & 0.24183 & 3.4355e-07 \\
    \cline{2-10}
        & Para. $y_0$                     & 0.06089 & 5.5206e-07 & 0.12065 & 4.3699e-06 & 0.18207 & 3.0924e-06 & 0.23964 & 1.2486e-06 \\
    \cline{2-10}
        & Distance                      & 0.00071 & 1.4810e-07 & 0.00183 & 1.5946e-06 & 0.00097 & 9.1203e-07 & 0.00218 & 4.1711e-07 \\
    \hline
    \multirow{3}{*}{$n=50$} & Inte. $y_0$ & 0.06149 & 4.8995e-08 & 0.12251 & 1.1746e-07 & 0.18300 & 4.2698e-07 & 0.24273 & 9.6758e-08 \\
    \cline{2-10}
        & Para. $y_0$                     & 0.06109 & 7.2222e-07 & 0.12152 & 6.1477e-07 & 0.18203 & 1.6656e-06 & 0.24172 & 3.8806e-07 \\
    \cline{2-10}
        & Distance                      & 0.00069 & 1.3734e-07 & 0.00099 & 2.3589e-07 & 0.00097 & 4.8985e-07 & 0.00101 & 1.3610e-07 \\
    \hline
  \end{tabular}
  }
  \label{tab:UtilityModel}
\end{table}

\begin{figure}[htbp]
    \centering
    \includegraphics[scale=0.45]{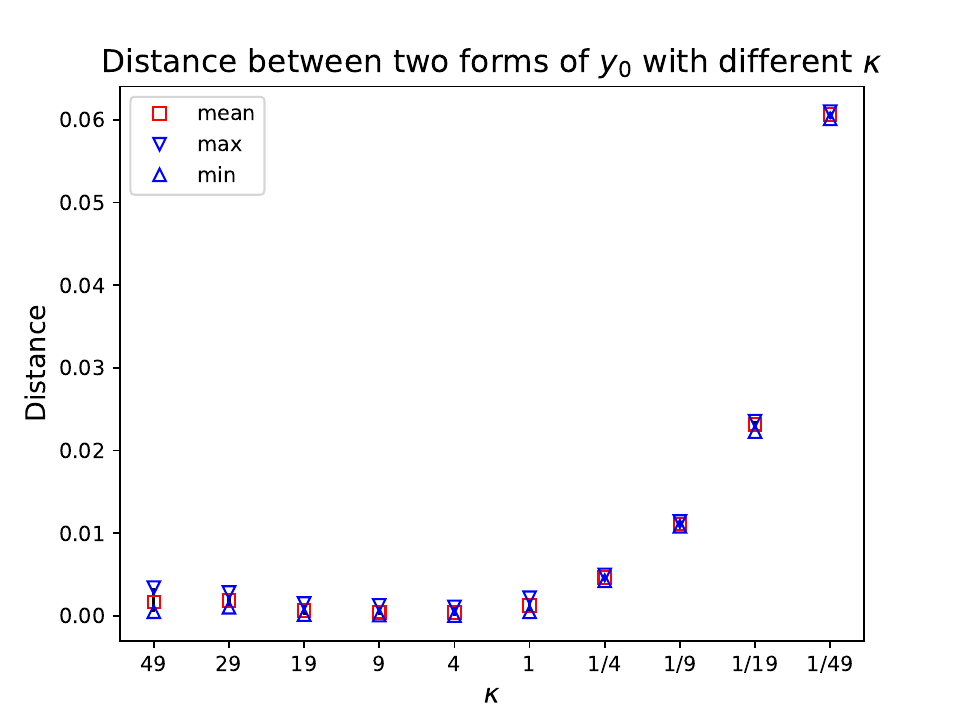}
    \includegraphics[scale=0.45]{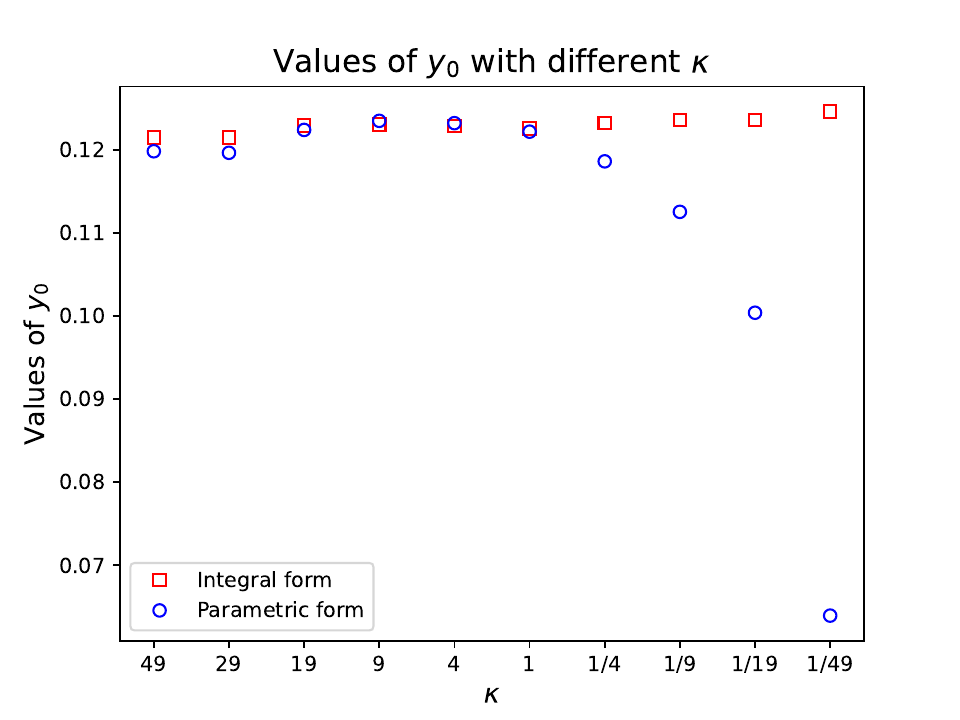}
    \caption{The distances of $y_0$ between the two different forms (the integral form and the parametric form) and the values of $y_0$ with different $\kappa$ for $n=10$ and $T=0.50$.}
    \label{fig-UtilityModel-k}
\end{figure}

We also vary the value of $\kappa$ for $n=10$ and $T=0.5$. The computation results are shown in \autoref{tab-Utility-k} and \autoref{fig-UtilityModel-k}, which also exhibit that the distances between the two different forms of $y_0$ are smaller and more stable when $\kappa > 1$.

\begin{table}[H]
  \centering
  \caption{The implementation results of $y_0$ with different $\kappa$ for case 2}
  \begin{tabular}{|c|c|c|c|c|c|c|}
    \hline
    \multicolumn{1}{|c|}{}& \multicolumn{2}{|c|}{Inte. $y_0$} & \multicolumn{2}{|c|}{Para. $y_0$} & \multicolumn{2}{|c|}{Distance} \\
    \hline
    & Mean & Var. & Mean & Var. & Mean & Var. \\
    \hline
    $\kappa=49$   & 0.12150 & 3.8599e-07 & 0.11983 & 2.0846e-06 & 0.00167 & 9.6376e-07 \\
    \hline
    $\kappa=29$   & 0.12150 & 8.3208e-08 & 0.11965 & 7.7199e-07 & 0.00185 & 3.5862e-07 \\
    \hline
    $\kappa=19$   & 0.12293 & 7.4394e-08 & 0.12240 & 6.7248e-07 & 0.00064 & 2.2293e-07 \\
    \hline
    $\kappa=9$    & 0.12307 & 1.8615e-07 & 0.12350 & 6.5421e-10 & 0.00043 & 2.0767e-07 \\
    \hline
    $\kappa=4$    & 0.12289 & 1.8920e-07 & 0.12322 & 1.3545e-09 & 0.00040 & 1.5793e-07 \\
    \hline
    $\kappa=1$    & 0.12255 & 1.7449e-06 & 0.12219 & 1.3545e-09 & 0.00124 & 1.5793e-07 \\
    \hline
    $\kappa=1/4$  & 0.12324 & 5.4207e-08 & 0.11862 & 2.6561e-08 & 0.00462 & 6.7300e-08 \\
    \hline
    $\kappa=1/9$  & 0.12362 & 1.2905e-09 & 0.11254 & 2.7437e-08 & 0.01108 & 3.5992e-08 \\
    \hline
    $\kappa=1/19$ & 0.12357 & 2.6459e-07 & 0.10039 & 9.0642e-09 & 0.02318 & 2.2544e-07 \\
    \hline
    $\kappa=1/49$ & 0.12461 & 1.0065e-07 & 0.06392 & 7.9881e-10 & 0.06069 & 1.1117e-07 \\
    \hline
  \end{tabular}
  \label{tab-Utility-k}
\end{table}

The above results demonstrate that our proposed method is effective for solving the stochastic recursive utility problem. Besides, we should do more optimizations on the follower's cost functionals, as verified by the implementation results.

\section{Conclusion}\label{sec:conclusion}

In this paper,
we propose a deep learning method for solving the stochastic control problem driven by fully-coupled FBSDEs.
We transform it into a stochastic Stackelberg differential game problem and propose a bi-level optimization method to solve this new problem.
The novel method has high flexibility and we can set the training parameters of the neural networks separately to optimize the leader's cost functional and the follower's cost functional independently.
In order to show the performance of our algorithm,
we give two optimal investment-consumption portfolio examples solved through the stochastic recursive utility models. The numerical results demonstrate remarkable performance.

\bibliographystyle{ieeetr}
\bibliography{ref}

\begin{thebibliography}{10}

\bibitem{bismut1973conjugate}
J.~Bismut, ``Conjugate convex functions in optimal stochastic control,'' {\em
  Journal of Mathematical Analysis and Applications}, vol.~44, no.~2,
  pp.~384--404, 1973.

\bibitem{pardoux1990adapted}
E.~Pardoux and S.~Peng, ``Adapted solution of a backward stochastic
  differential equation,'' {\em Systems \& Control Letters}, vol.~14, no.~1,
  pp.~55--61, 1990.

\bibitem{peng1997financial}
N.~El~Karoui, S.~Peng, and M.~C. Quenez, ``Backward stochastic differential
  equations in finance,'' {\em Mathematical finance}, vol.~7, no.~1, pp.~1--71,
  1997.

\bibitem{peng1990stochasticmax}
S.~Peng, ``A general stochastic maximum principle for optimal control
  problems,'' {\em Siam Journal on Control and Optimization}, vol.~28, no.~4,
  pp.~966--979, 1990.

\bibitem{wu1999Fully}
S.~Peng and Z.~Wu, ``Fully coupled forward-backward stochastic differential
  equations and applications to optimal control,'' {\em Siam Journal on Control
  and Optimization}, vol.~37, no.~3, pp.~825--843, 1999.

\bibitem{antonelli1993backward}
F.~Antonelli, ``Backward-forward stochastic differential equations,'' {\em
  Annals of Applied Probability}, vol.~3, no.~3, pp.~777--793, 1993.

\bibitem{ma1996hedging}
J.~Cvitani{\'c} and J.~Ma, ``Hedging options for a large investor and
  forward-backward sde's,'' {\em The annals of applied probability}, vol.~6,
  no.~2, pp.~370--398, 1996.

\bibitem{hu1995solution}
Y.~Hu and S.~Peng, ``Solution of forward-backward stochastic differential
  equations,'' {\em Probability Theory and Related Fields}, vol.~103, no.~2,
  pp.~273--283, 1995.

\bibitem{pardoux1999forward}
E.~Pardoux and S.~Tang, ``Forward-backward stochastic differential equations
  and quasilinear parabolic pdes,'' {\em Probability Theory and Related
  Fields}, vol.~114, no.~2, pp.~123--150, 1999.

\bibitem{ma1994solving}
J.~Ma, P.~Protter, and J.~Yong, ``Solving forward-backward stochastic
  differential equations explicitly—a four step scheme,'' {\em Probability
  theory and related fields}, vol.~98, no.~3, pp.~339--359, 1994.

\bibitem{peng1993backward}
S.~Peng, ``Backward stochastic differential equations and applications to
  optimal control,'' {\em Applied Mathematics and Optimization}, vol.~27,
  no.~2, pp.~125--144, 1993.

\bibitem{dokuchaev1999JMAA}
N.~Dokuchaev and X.~Y. Zhou, ``Stochastic controls with terminal contingent
  conditions,'' {\em Journal of Mathematical Analysis and Applications},
  vol.~238, no.~1, pp.~143--165, 1999.

\bibitem{ji2006cis}
S.~Ji and X.~Y. Zhou, ``{A maximum principle for stochastic optimal control
  with terminal state constraints, and its applications},'' {\em Communications
  in Information and Systems}, vol.~6, no.~4, pp.~321 -- 338, 2006.

\bibitem{Yong_stochastic_control}
J.~Yong and X.~Zhou, {\em Stochastic Controls-Hamiltonian System and HJB
  Equations}.
\newblock Springer, 1999.

\bibitem{hu2017puqr}
M.~Hu, ``Stochastic global maximum principle for optimization with recursive
  utilities,'' {\em Probability, Uncertainty and Quantitative Risk}, vol.~2,
  no.~1, 2017.

\bibitem{hu2018global}
M.~Hu, S.~Ji, and X.~Xue, ``A global stochastic maximum principle for fully
  coupled forward-backward stochastic systems,'' {\em SIAM Journal on Control
  and Optimization}, vol.~56, no.~6, pp.~4309--4335, 2018.

\bibitem{li2014siam-jco}
J.~Li and Q.~Wei, ``Optimal control problems of fully coupled fbsdes and
  viscosity solutions of hamilton--jacobi--bellman equations,'' {\em SIAM
  Journal on Control and Optimization}, vol.~52, no.~3, pp.~1622--1662, 2014.

\bibitem{hu2019siamCO}
M.~Hu, S.~Ji, and X.~Xue, ``The existence and uniqueness of viscosity solution
  to a kind of hamilton--jacobi--bellman equation,'' {\em SIAM Journal on
  Control and Optimization}, vol.~57, no.~6, pp.~3911--3938, 2019.

\bibitem{hu2020esaim-cocv}
M.~Hu, S.~Ji, and X.~Xue, ``Stochastic maximum principle, dynamic programming
  principle, and their relationship for fully coupled forward-backward
  stochastic controlled systems,'' {\em ESAIM: COCV}, vol.~26, p.~81, 2020.

\bibitem{han2016deep}
J.~Han and W.~E, ``Deep learning approximation for stochastic control
  problems,'' {\em NIPS Workshop on Deep Reinforcement Learning}, 2016.

\bibitem{carmona2021convergence2}
R.~Carmona and M.~Laurière, ``Convergence analysis of machine learning
  algorithms for the numerical solution of mean field control and games: ii --
  the finite horizon case,'' {\em arXiv preprint arXiv:1908.01613}.

\bibitem{numericalPDE2012review}
E.~Tadmor, ``A review of numerical methods for nonlinear partial differential
  equations,'' {\em Bulletin of the American Mathematical Society}, vol.~49,
  no.~4, pp.~507--554, 2012.

\bibitem{regression_for_BSDEs}
B.~Bouchard and N.~Touzi, ``Discrete-time approximation and monte-carlo
  simulation of backward stochastic differential equations,'' {\em Stochastic
  Processes \& Their Applications}, vol.~111, no.~2, pp.~175--206, 2004.

\bibitem{FBSDE_Ma2008}
M.~A. Jin, S.~Jie, and Y.~Zhao, ``On numerical approximations of
  forward-backward stochastic differential equations,'' {\em Siam Journal on
  Numerical Analysis}, vol.~46, no.~5, pp.~2636--2661, 2008.

\bibitem{time_discretize_FBSDE_Zhang}
C.~Bender and J.~Zhang, ``Time discretization and markovian iteration for
  coupled fbsdes,'' {\em Annals of Applied Probability}, vol.~18, no.~1,
  pp.~143--177, 2008.

\bibitem{Zhao2016Multistep}
Y.~Fu, W.~Zhao, and T.~Zhou, ``Multistep schemes for forward backward
  stochastic differential equations with jumps,'' {\em Journal of Scientific
  Computing}, vol.~69, no.~2, pp.~1--22, 2016.

\bibitem{Milstein2004Numerical}
G.~N. Milstein and M.~V. Tretyakov, ``Numerical algorithms for forward-backward
  stochastic differential equations connected with semilinear parabolic
  equations,'' vol.~28, pp.~561--582, 2004.

\bibitem{Yu2016Efficient}
Y.~Fu, W.~Zhao, and Z.~Tao, ``Efficient spectral sparse grid approximations for
  solving multi-dimensional forward backward sdes,'' {\em Discrete and
  Continuous Dynamical Systems - Series B}, vol.~22, no.~9, 2016.

\bibitem{huijskens2016efficient}
T.~P. Huijskens, M.~Ruijter, and C.~W. Oosterlee, ``Efficient numerical fourier
  methods for coupled forward-backward sdes,'' {\em Journal of Computational
  and Applied Mathematics}, vol.~296, pp.~593--612, 2016.

\bibitem{WeinanDLforHigh_dim}
J.~Han, A.~Jentzen, and W.~E, ``Solving high-dimensional partial differential
  equations using deep learning,'' {\em Proceedings of the National Academy of
  Sciences of the United States of America}, vol.~115, no.~34, pp.~8505--8510,
  2018.

\bibitem{WeinanDLforBSDE}
W.~E, J.~Han, and A.~Jentzen, ``Deep learning-based numerical methods for
  high-dimensional parabolic partial differential equations and backward
  stochastic differential equations,'' {\em Communications in Mathematics \&
  Statistics}, vol.~5, no.~4, pp.~349--380, 2017.

\bibitem{deeplearning_FBSDE}
J.~Han and J.~Long, ``Convergence of the deep bsde method for coupled fbsdes,''
  {\em arXiv:1811.01165}, 2018.

\bibitem{Peng_FBSDE_numerical}
S.~Ji, S.~Peng, Y.~Peng, and X.~Zhang, ``Three algorithms for solving
  high-dimensional fully coupled fbsdes through deep learning,'' {\em IEEE
  Intelligent Systems}, vol.~35, no.~3, pp.~71--84, 2020.

\bibitem{hure_deep_2020}
C.~Huré, H.~Pham, and X.~Warin, ``Deep backward schemes for high-dimensional
  nonlinear {PDEs},'' {\em arXiv:1902.01599 [cs, math, stat]}, June 2020.

\bibitem{pham_neural_2021}
H.~Pham, X.~Warin, and M.~Germain, ``Neural networks-based backward scheme for
  fully nonlinear {PDEs},'' {\em SN Partial Differential Equations and
  Applications}, vol.~2, p.~16, 2021.

\bibitem{beck_machine_2019}
C.~Beck, W.~E, and A.~Jentzen, ``Machine {Learning} {Approximation}
  {Algorithms} for {High}-{Dimensional} {Fully} {Nonlinear} {Partial}
  {Differential} {Equations} and {Second}-order {Backward} {Stochastic}
  {Differential} {Equations},'' {\em Journal of Nonlinear Science}, vol.~29,
  no.~4, pp.~1563--1619, 2019.

\bibitem{duffie1992stochastic}
D.~Duffie and L.~G. Epstein, ``Stochastic differential utility,'' {\em
  Econometrica: Journal of the Econometric Society}, pp.~353--394, 1992.

\bibitem{el2001dynamic}
N.~El~Karoui, S.~Peng, and M.~C. Quenez, ``A dynamic maximum principle for the
  optimization of recursive utilities under constraints,'' {\em Annals of
  applied probability}, pp.~664--693, 2001.

\bibitem{chen2002ambiguity}
Z.~Chen and L.~Epstein, ``Ambiguity, risk, and asset returns in continuous
  time,'' {\em Econometrica}, vol.~70, no.~4, pp.~1403--1443, 2002.

\bibitem{buckdahn2010existence}
R.~Buckdahn, B.~Labed, C.~Rainer, and L.~Tamer, ``Existence of an optimal
  control for stochastic control systems with nonlinear cost functional,'' {\em
  Stochastics: An International Journal of Probability and Stochastics
  Processes}, vol.~82, no.~3, pp.~241--256, 2010.

\bibitem{bahlali2011existence}
K.~Bahlali, B.~Gherbal, and B.~Mezerdi, ``Existence of optimal controls for
  systems driven by fbsdes,'' {\em Systems \& control letters}, vol.~60, no.~5,
  pp.~344--349, 2011.

\bibitem{bahlali2018existence}
K.~Bahlali, O.~Kebiri, B.~Mezerdi, and A.~Mtiraoui, ``Existence of an optimal
  control for a coupled fbsde with a non degenerate diffusion coefficient,''
  {\em Stochastics}, vol.~90, no.~6, pp.~861--875, 2018.

\bibitem{KohlmannZhou2000}
M.~Kohlmann and X.~Zhou, ``Relationship between backward stochastic
  differential equations and stochastic controls: a linear-quadratic
  approach,'' {\em SIAM Journal on Control and Optimization}, vol.~38, no.~5,
  pp.~1392--1392, 2000.

\bibitem{LimZhou2001}
A.~Lim and X.~Zhou, ``Linear-quadratic control of backward stochastic
  differential equations,'' {\em SIAM Journal on Control and Optimization},
  2001.

\bibitem{yong2002leader}
J.~Yong, ``A leader-follower stochastic linear quadratic differential game,''
  {\em SIAM Journal on Control and Optimization}, vol.~41, no.~4,
  pp.~1015--1041, 2002.

\bibitem{pham2018deep_2}
A.~Bachouch, C.~Hure, N.~Langrene, and H.~Pham, ``Deep neural networks
  algorithms for stochastic control problems on finite horizon, part 2:
  numerical applications,'' {\em arXiv: 1812.05916v1}, 2018.

\bibitem{ji2021solving}
S.~Ji, S.~Peng, Y.~Peng, and X.~Zhang, ``Solving stochastic optimal control
  problem via stochastic maximum principle with deep learning method,'' 2021.

\end{thebibliography}

\end{document}